\documentclass[twocolumn]{autart}    
\usepackage{indentfirst} 
\usepackage{graphicx}     
\usepackage{amssymb}
\usepackage{CJK}
\usepackage{bm}
\usepackage{subfigure}
\usepackage{multirow}
\usepackage{booktabs}
\usepackage{pifont}
\usepackage{textcomp}
\usepackage{algorithm}
\usepackage{caption}
\usepackage{algpseudocode}
\usepackage[dvips]{epsfig}    
\usepackage{amsmath}

\newcounter{counter}

\begin{document}
\begin{frontmatter}

\title{Local Differential Privacy via Dynamic Quantization in Distributed Online Stochastic Optimization\thanksref{footnoteinfo}} 

\thanks[footnoteinfo]{The material in this paper was not presented at any conference. Corresponding author: Qian Ma.}

\author[NJUST]{Zhiguo Zhang}\ead{zgzhang0426@163.com},    
\author[NJUST]{Cheng Kui}\ead{kuicheng404@gmail.com}, 
\author[NJUST]{Qian Ma}\ead{qma@njust.edu.cn},
\author[HUST]{Dongrui Wu}\ead{drwu09@gmail.com}

\address[NJUST]{School of Automation, Nanjing University of Science and Technology, Nanjing 210094, China}  
\address[HUST]{Key Laboratory of the Ministry of Education for 
Image Processing and Intelligent Control, School of Artificial Intelligence and Automation, Huazhong University of Science and Technology, Wuhan 430074, China}


\begin{abstract}                          
Distributed online stochastic optimization has received extensive attention in large-scale distributed learning and other related fields due to its unique advantage in processing streaming data.
However, information exchange through the communication network during the optimization process may lead to privacy leakage.
To address this issue, this paper proposes a locally differentially private distributed online stochastic optimization algorithm that employs an elaborately designed dynamic stochastic quantizer to mask the exchanged information prior to communication.
Theoretical analysis shows that the proposed algorithm not only converges almost surely to the optimal solution but also achieves $(0,\delta^i)$-local differential privacy for each agent $i$ even when the number of iterations tends to infinity.
Furthermore, the algorithm is fully distributed and applicable to scenarios where the interaction network among agents is a directed graph.
To the best of our knowledge, this is the first work on distributed online stochastic optimization that simultaneously achieves exact convergence and rigorous local differential privacy over a directed graph by exploiting quantization effects.
Numerical experiments of distributed online training on the mushroom classification dataset, handwritten digits recognition dataset, and brain-computer interface dataset verify the effectiveness of the proposed method.
\end{abstract}

\end{frontmatter}

\section{Introduction}
Distributed optimization has received considerable attention due to its significant role in smart grids, sensor networks, large-scale machine learning, and other areas \cite{Mohiuddin_Optimal,Predd_Distributed,Qian_Distributed,Zhao_Scalable,Zhu_On}.
In a typical distributed optimization problem, each agent cooperates with its neighbors by exchanging information to minimize (or maxmize) the sum of all local cost (or payoff) functions, where each agent has access only to its own local  function.
Note that precise gradients are often unavailable or computationally expensive in many practical applications, which has driven the development of distributed stochastic optimization.
Many significant works have been proposed, see \cite{Lu_Convergence,Lei_Distributed,Huang_Distributed,Zhao_Confidence,Yi_Zeroth} and the references therein.
Among them, reformulating distributed stochastic optimization as an empirical risk minimization (ERM) problem allows for the approximate estimation of the original optimal solution from observed data.
However, the traditional ERM method requires all data to be available in advance, which makes it incapable of handling sequentially arriving data.
This gave rise to distributed online stochastic optimization, and numerous achievements have been made in this field, such as \cite{Lee_Stochastic,Zhou_Online,Cao_Online,Xu_Online,Yang_Online}.
\\
\indent\setlength{\parindent}{1em}
Notably, the information exchange process over the communication network in distributed optimization may lead to privacy leakage.
For example, while multiple power plants in a smart grid coordinate their generation capacities to minimize total electricity costs, they need to protect the privacy of their generating information, which is sensitive in energy sales bidding.
In fact, attackers are able to recover the exact original data from the shared gradients or model updates in the absence of privacy-preserving mechanisms \cite{Zhu_Deep}.
Therefore, it is of great significance to develop distributed optimization algorithms with privacy protection.
A homomorphic encryption-based distributed projection gradient algorithm was proposed in \cite{Lu_Privacy}, allowing participants to keep the shared information encrypted throughout the optimization procedure. Neverthless, the homomorphic encryption technique involves complex ciphertext computations, leading to high computational and communication costs.
In \cite{Lou_Privacy}, the privacy-preserving distributed optimization was achieved by adding a constant uncertain parameter to the projection step.
In \cite{Gade_Private}, correlated noise was injected to obfuscate the information transmitted between agents, which reduces the risk of privacy leakage in distributed optimization.
However, since the randomness introduced in \cite{Lou_Privacy} and \cite{Gade_Private} is interrelated, the level of the privacy protection is limited.
\\
\indent\setlength{\parindent}{1em}
Recently, differential privacy (DP) has gradually become the mainstream approach for privacy protection in distributed optimization due to its strong resistance to arbitrary post-processing and auxiliary information \cite{Dwork_The}.
The common implementation of DP involves injecting irrelevant noise into the exchanged information to prevent eavesdroppers from accurately inferring the private information of individual agents.
While significant progress has been made in differentially private distributed optimization (see \cite{Liu_Distributed,Huang_Differential,Lü_Privacy-preserving,Gratton_Privacy-preserved}), a common limitation of these works is the trade-off between convergence accuracy and privacy protection level. In other words, strong privacy protection is achieved at the expense of convergence accuracy.
By introducing a weakening factor into iteration to mitigate the impact of DP noise, Wang et al. in \cite{Wang_Tailoring} simultaneously achieved accurate convergence and a finite privacy budget, even with an infinite number of iterations. Subsequently, many other works have emerged based on this approach, as seen in \cite{Wang_Robust,Chen_Differentially,Chen_Locally}.
More importantly, \cite{Chen_Locally} adopted a local differential privacy (LDP) framework, which is a local (distributed) model of DP designed for scenarios where all communication channels may be compromised by malicious attackers and no agent is trustworthy.
Compared with the traditional DP framework in \cite{Liu_Distributed,Huang_Differential,Lü_Privacy-preserving,Gratton_Privacy-preserved,Wang_Tailoring,Wang_Robust,Chen_Differentially}, 
this LDP framework provides a significantly more stringent privacy guarantee under a much weaker trust model.
\\
\indent\setlength{\parindent}{1em}
In addition to adding random noise, quantization effects offer another technical approach that not only achieves DP but also reduces communication overhead.
In the pioneering work \cite{Wang_Quantization}, Wang et al. proposed a distributed stochastic optimization algorithm that achieves DP via quantization.
Specifically, by introducing a random ternary quantizer to quantize any exchanged information into three numerical levels before transmission, precise convergence and strict $(0,\delta)$-DP guarantee at each iteration were achieved by the proposed algorithm. However,
it is worth noting that this algorithm will degenerate to $(0,1)$-DP as the number of iterations tends to infinity, which is equivalent to directly outputting sensitive information.
Moreover, the quantizer's stepsize in \cite{Wang_Quantization} is required to be larger than the infinite norm of the decision variable, which is difficult to guarantee in practice.
To address the privacy leakage in cooperative bandit games, a scheme combining stochastic quantization with a binary stochastic response was adopted in \cite{Lin_Quantization} to mask the action estimation before communication. This scheme achieves the same order of regret bound as the non-private setting while maintaining strict $(\epsilon,0)$-DP at each iteration.
However, the aforementioned algorithm cannot ensure a finite privacy budget over infinite iterations, meaning that sensitive information may no longer be effectively protected as the iteration count grows.
Furthermore, both \cite{Wang_Quantization} and \cite{Lin_Quantization} adopt the traditional DP framework, which implicitly assumes mutual trust among agents and the presence of a centralized data aggregator capable of determining the amount of injected noise.
This makes it inapplicable to fully distributed and fully untrusted environments.
In light of the above analysis, how to leverage quantization effects to achieve both precise convergence in distributed online stochastic optimization and LDP for each agent over infinite iterations remains an open and challenging issue.
\\
\indent\setlength{\parindent}{1em}
In this paper, the LDP-based distributed online stochastic optimization over directed graph  is discussed.
A stochastic quantizer with a dynamically decaying stepsize is designed to obfuscate the information exchanged between agents, thereby achieving privacy protection.
By employing gradient tracking and online eigenvector estimation techniques, we propose a locally differentially private fully distributed online stochastic optimization algorithm over unbalanced directed graph, and prove its almost sure convergence to the optimal solution and rigorous $(0,\delta^i)$-LDP for each agent $i$.
The main contributions are as follows:
\\
\indent\setlength{\parindent}{1em}
1) Unlike the static quantizers used in \cite{Wang_Quantization} and \cite{Lin_Quantization}, a dynamic stochastic quantizer whose quantization stepsize gradually decreases as the number of iterations increases is developed. Instead of relying on the weakening factor to reduce the intensity of information exchange between agents to ensure accurate convergence, we leverage the quantizer's inherent variance decay to avoid the tradeoff between privacy strength and convergence accuracy, without compromising convergence speed.
\\
\indent\setlength{\parindent}{1em}
2) Compared with \cite{Wang_Quantization} and \cite{Lin_Quantization} that can only achieve DP at each iteration, our distributed online stochastic optimization algorithm based on the proposed dynamic stochastic quantizer can achieve $(0,\delta^i)$-LDP over an infinite number of iterations while maintaining the convergence accuracy. Even when the problem degenerates to the standard DP setting, our approach still ensures $(0,\delta)$-DP for an infinite time horizon, which is strictly stronger than the per-iteration privacy guarantees of \cite{Wang_Quantization} and \cite{Lin_Quantization}.
\\
\indent\setlength{\parindent}{1em}
3) In contrast to existing works on differentially private distributed optimization \cite{Liu_Distributed,Huang_Differential,Lü_Privacy-preserving,Gratton_Privacy-preserved,Wang_Tailoring,Wang_Robust,Chen_Differentially,Chen_Locally}, the proposed privacy-preserving framework is more flexible and practical. On the one hand, the adopted LDP framework allows each agent to freely choose its required privacy strength, and avoids the reliance on any trusted third parties. On the other hand, the introduced quantizer can save network communication resources, which is particularly important in scenarios with limited communication bandwidth. 

\noindent\textbf{Notations.} We use $\mathbb{R}^d$ to denote the Euclidean space of dimension $d$. Let $\otimes$ denote the Kronecker product. We write $I_d$ for the identity matrix of dimension $d$, $\mathbf{0}_d$ and $\mathbf{1}_d$ for the $d\text{-}$dimensional column vectors with all entries equal to $0$ and $1$, respectively.
$\left\langle \cdot,\cdot \right\rangle$ denotes the inner product.
For a given vector $a$, $\text{diag}(a)$ represents a diagonal matrix whose diagonal entries are the components of the vector $a$ in sequence, i.e., $\text{diag}(a)=\text{diag}\{a_1,a_2,\cdots,a_n\}$.
For a given matrix $P$, we denote its Euclidean norm, inverse and transpose as $\left\|P\right\|$, $P^{-1}$, and $P^T$, respectively. We abbreviate almost surely by a.s., and independent and identically distributed by i.i.d..

\section{Problem Formulation and Preliminaries}
\subsection{Distributed Online Stochastic Optimization}
In this paper, we consider the following distributed stochastic optimization problem
\begin{equation}
  \min_{\theta\in\mathbb{R}^d}F(\theta)=\frac{1}{m}\sum_{i=1}^{m}f_i(\theta),\ f_i(\theta)=\mathbb{E}_{x^i\sim\Phi ^i}[l_i(\theta,x^i)], \label{problem1}
\end{equation}
where the local objective function $f_i(\theta): \mathbb{R}^d \to \mathbb{R}$ represents the mathematical expectation of local loss function $l_i(\theta,x^i)$, and $x^i$ is a random variable drawn from an unknown probability distribution $\Phi^i$.
To accommodate online scenarios, we consider a network of $m$ agents that cooperatively learn an optimal solution $\theta^*$ to problem \eqref{problem1} with data acquired in serial, like \cite{Chen_Local}. To this end, we reformulate \eqref{problem1} as the following form
\begin{equation}
  \min_{\theta \in \mathbb{R}^d} F_t(\theta) = \frac{1}{m} \sum_{i=1}^{m} f_t^i(\theta),\ f_t^i(\theta) = \frac{\sum_{k=0}^{t} l_i(\theta, x_k^i)}{t+1},\label{problem2}
\end{equation}
where $x_k^i$ represents the data-point acquired by agent $i$ at time $k$.
To facilitate subsequent analysis, the following standard assumptions are made.

\setcounter{counter}{0}
\refstepcounter{counter}
\label{assumption 1}
\noindent\textbf{Assumption 1.}
Problem \eqref{problem1} has an optimal solution $\theta^*$. The objective function $F(\cdot)$ is convex and has bounded level sets.

\setcounter{counter}{1}
\refstepcounter{counter}
\label{assumption 2}
\noindent\textbf{Assumption 2.}
For any $T \in \mathbb{N}^+$, the random data-points $\{ x_k^i\}$ of agent $i$ are i.i.d. across different time instants $k=0,\cdots,T-1$. In addition, the following statements hold: (i) $\mathbb{E}[\nabla l_i(\theta, x_k^i)]=\nabla f_i(\theta)$; (ii) $\mathbb{E}[\| \nabla l_i(\theta, x_k^i) - \nabla f_i(\theta)\|^2] \leq \kappa^2$; and (iii) there exists a Lipschitz constant $L>0$ such that $\| \nabla l_i(\theta_1,x_k^i)-\nabla l_i(\theta_2,x_k^i)\| \leq L\|\theta_1 - \theta_2\|$ for any $\theta_1,\theta_2 \in \mathbb{R}^d$.

\setcounter{counter}{2}
\refstepcounter{counter}
\label{assumption 3}
\noindent\textbf{Assumption 3.}
For any $T \in \mathbb{N}^+$, there exists some positive constant $d_l$ such that $\|\nabla l_i(\theta,x_k^i)\|_1 \leq d_l$ holds for any $\theta\in\mathbb{R}^d$,  $i\in[m]$ and $k=0,\cdots,T-1$.

These assumptions are standard for convergence analysis in distributed stochastic optimization (see, e.g., \cite{Wang_Quantization}, \cite{Pu_A}, and \cite{Bullins_A}).
Then, an important definition and some auxiliary lemmas are provided.

\setcounter{counter}{0}
\refstepcounter{counter}
\label{definition 1}
\noindent\textbf{Definition 1 \cite{Pu_Push-pull}.}
Given an arbitrary vector norm $\|\cdot\|_A$ on $\mathbb{R}^d$, for any $w \in \mathbb{R}^{m\times d}$, we define
\begin{equation*}
  \|w\|_A \triangleq \left\|(\left\|w(1)\right\|_A,\left\|w(2)\right\|_A, \cdots, \left\|w(m)\right\|_A)\right\|_2,
\end{equation*}
where $w(1),w(2),\cdots,w(m)\in\mathbb{R}^d$ are columns of $w$, and the subscript $2$ denotes the $2$-norm.

\setcounter{counter}{0}
\refstepcounter{counter}
\label{lemma 1}
\noindent\textbf{Lemma 1 \cite{Wang_Tailoring}.}
Let $\{ v_t \}$, $\{ \alpha_t \}$, and $\{ p_t \}$ be random nonnegative sequences, and $\{ q_t \}$ be a deterministic nonnegative scalar sequence satisfying  $\sum_{t=0}^{\infty}\alpha_t < \infty$ a.s., $\sum_{t=0}^{\infty}q_t = \infty$, $\sum_{t=0}^{\infty}p_t < \infty$ a.s., and the following inequality
\begin{eqnarray*}
  \mathbb{E}[v_{t+1}|\mathcal{F}_t] \leq (1 + \alpha_t - q_t)v_t + p_t,\ \forall t \ge 0 \ \ \text{a.s.},
\end{eqnarray*}
where $\mathcal{F}_t = \{ v_\ell, \alpha_\ell, p_\ell; 0 \leq \ell \leq t\}$. Then, $\sum_{t=0}^{\infty}q_t v_t \leq \infty$ and $\lim_{t \to \infty}v_t = 0$ hold a.s..

\setcounter{counter}{1}
\refstepcounter{counter}  
\label{lemma 2}
\noindent\textbf{Lemma 2 \cite{Wang_Tailoring}.}
Let $\{ \boldsymbol{v_t}\} \subset \mathbb{R}^d$ and $\{ \boldsymbol{u_t}\} \subset \mathbb{R}^p$ be random nonnegative vector sequences, $\{ a_t \}$ and $\{ b_t \}$ be random nonnegative scalar sequences, and $\{ V_t\}$ and $\{ H_t\}$ be random sequences of nonnegative matrices such that
\begin{eqnarray*}
  \mathbb{E}[\boldsymbol{v_{t+1}} | \mathcal{F}_t] \leq (V_t + a_t \mathbf{11}^T)\boldsymbol{v_t} + b_t \mathbf{1} - H_t \boldsymbol{u_t},\ \forall t \ge 0\ \ \text{a.s.},
\end{eqnarray*}
where $\mathcal{F}_t = \{ \boldsymbol{v_\ell}, \boldsymbol{u_\ell}, a_\ell, b_\ell, V_\ell, H_\ell; 0 \leq \ell \leq t\}$. Assume that $\{ a_t\}$ and $\{ b_t\}$ satisfy $\sum_{t=0}^{\infty}a_t < \infty$ and $\sum_{t=0}^{\infty}b_t < \infty$ a.s., and that there exists a vector $\pi > 0$ such that $\pi^T V_t \leq \pi^T$ and $\pi^T H_t \ge 0 $ hold a.s. for all $t \ge 0$. Then, we have (i) $\{\pi^T \boldsymbol{v_t}\}$ converges to some random variable $\pi^T \boldsymbol{v} \ge 0$ a.s.; (ii) $\{\boldsymbol{v_t}\}$ is bounded a.s.; (iii) $\sum_{t=0}^{\infty}\pi^T H_t \boldsymbol{u_t} < \infty$ holds a.s..

\subsection{Local Differential Privacy}
As shown in \cite{Wang_Quantization}, two adversary models are commonly considered in the privacy protection problem of distributed optimization:
\begin{enumerate}
  \item [$\bullet$] \emph{A honest-but-curious adversary}. This adversary is defined as an agent that follows the prescribed protocol and correctly computes its iterative states while having access to certain internal states from others. However, its objective is to infer the sensitive information of other agents.
  \item [$\bullet$] \emph{An eavesdropper}. This adversary refers to an external entity capable of eavesdropping on and monitoring all communication channels, thereby capturing messages transmitted by any agent. This allows the eavesdropper to infer sensitive information of internal agents.
\end{enumerate}

In order to address the privacy leakage issues caused by adversaries, LDP is introduced as a privacy notion that requires each agent to randomize its own messages before sharing. As a result, privacy is protected even when no trusted data curator exists and all communication channels are compromised.
Before providing the definition of LDP, we first introduce the concept of adjacency on the local dataset of agent $i$.

\setcounter{counter}{1}
\refstepcounter{counter}  
\label{definition 2}
\noindent\textbf{Definition 2 (Adjacency) \cite{Chen_Locally}.}
For any $T \in \mathbb{N}^+$ and any agent $i \in [m]$, given two local datasets $\mathcal{D}^i = \{x_0^i, \cdots, x_{T-1}^i\}$ and $\mathcal{D}'^{i} = \{x_0'^{i}, \cdots, x_{T-1}'^{i}\}$, $\mathcal{D}^i$ is said to be adjacent to $\mathcal{D}'^{i}$ if there exists a time instant $k=0,\cdots,T-1$ such that $x_k^i \ne x_k'^{i}$ while $x_p^i = x_p'^{i}$ for all $p=0,\cdots,T-1$ and $p \ne k$.

According to Definition \ref{definition 2}, two local datasets $\mathcal{D}^i$ and $\mathcal{D}'^i$ are adjacent if and only if they differ by exactly one element.
Following \cite{Li_Convergence} and denoting this adjacency relationship between $\mathcal{D}^i$ and $\mathcal{D}'^{i}$ as $\mathrm{Adj}(\mathcal{D}^i, \mathcal{D}'^{i})$, we present the definition of $(0,\delta^i)$-LDP for each agent $i$ as follows.

\setcounter{counter}{2}
\refstepcounter{counter}  
\label{definition 3}
\noindent\textbf{Definition 3 ($(0, \delta^i)$-LDP).}
Given $0 \leq \delta^i \leq 1$, a mechanism $\mathcal{M}$ achieves $(0, \delta^i)$-LDP for $\mathrm{Adj}(\mathcal{D}^i, \mathcal{D}'^{i})$ if $\mathbb{P}[\mathcal{M}(\mathcal{D}^i) \in \mathcal{S}] \leq \mathbb{P}[\mathcal{M}(\mathcal{D}'^{i}) \in \mathcal{S}] + \delta^i$ for any Borel-measurable set $\mathcal{S} \subseteq \mathrm{Range}(\mathcal{M})$.

The above definition of $(0,\delta^i)$-LDP ensures that an adversary with access to all exchanged messages in the network cannot infer any agent's local dataset with significant probability.
A smaller $\delta^i$ corresponds to a stronger level of privacy protection.
It is worth noting that the considered notion of $(0,\delta^i)$-LDP is more stringent than $(\epsilon^i,\delta^i)$-LDP \cite{Li_Convergence} in the local model.

\subsection{Dynamic Stochastic Quantizer}
Due to the privacy issue caused by the information exchange among agents, a dynamic stochastic quantizer $Q$ is presented to quantize the exchanged information in this paper. For any finite input $y \in \mathbb{R}$, we can express it as $y = z + nd_t$, where $d_t$ is the quantization stepsize decaying as $t$ increases, and $n$ is the largest integer such that $z \in (0,d_t]$. The quantizer $Q$ randomly maps the input $y$ according to the following probabilities
\begin{equation}
  Q(y)=
  \begin{cases}
  nd_t & \mbox{with probability}\ 1 - \frac{z}{d_t}, \\
  (n+1)d_t & \mbox{with probability}\quad \frac{z}{d_t}.
  \end{cases}
  \label{quantizer}
\end{equation}
It can be easily calculated that $\mathbb{E}[Q(y)|y]=y$, and $\mathbb{E}[(y-Q(y))^2|y] \leq \frac{d_t^2}{4}$.

\section{Locally Differentially Private Distributed Online Stochastic Optimization Algorithm}
In this section, we propose a fully distributed algorithm for locally differentially private online stochastic optimization over unbalanced directed graph and prove that it convergences almost surely to the optimal solution despite the presence of quantization errors, and each agent $i\in [m]$ can achieve rigorous $(0, \delta^i)$-LDP.
In our algorithm, each agent $i$ maintains and updates three variables at iteration $t$, including an optimization variable $\theta_t^i$, a gradient-tracking variable $\psi_t^i$, and an eigenvector-estimation variable $z_t^i$. These variables will be shared with its neighbors by two different communication networks called $\mathcal{G}_R$ and $\mathcal{G}_C$, which are induced by matrices $R=[R_{ij}] \in \mathbb{R}^{m \times m}$ and $C=[C_{ij}] \in \mathbb{R}^{m \times m}$, respectively. Specifically, a directed link $(i,j)$ from agent $j$ to agent $i$ exists if and only if $R_{ij} > 0$; the same holds for $\mathcal{G}_C$ with respect to $C$. We make the following assumption on $R$ and $C$. Note that $\mathcal{G}_{C^T}$ is equivalent to $\mathcal{G}_C$ with edge directions reversed.

\setcounter{counter}{3}
\refstepcounter{counter}  
\label{assumption 4}
\noindent\textbf{Assumption 4.}
The matrices $R$ and $C$ have nonnegative off-diagonal entries, i.e., $R_{ij} \ge 0\text{ and }C_{ij} \ge 0\text{ for all }i \ne j$, and their diagonal entries are negative and satisfy $1+R_{ii}>0$ and $1+C_{ii}>0$ with $R_{ii} = - \sum_{j \in \mathbb{N}_{R,i}^\text{in}}R_{ij}$ and $C_{ii} = -\sum_{j \in \mathbb{N}_{C,i}^\text{out}}C_{ji}$, which implies that $R$ has zero row sums and $C$ has zero column sums. Moreover, the induced graphs $\mathcal{G}_R$ and $\mathcal{G}_{C^T}$ satisfy (i) $\mathcal{G}_R$ and $\mathcal{G}_{C^T}$ each contain at least one spanning tree;
(ii) There exists at least one agent that serves as a common root for spanning trees in both $\mathcal{G}_R$ and $\mathcal{G}_{C^T}$.

\noindent\textbf{Remark 1.}
This assumption on $\mathcal{G}_R$ and $\mathcal{G}_{C^T}$ is weaker than the condition that the induced graphs of both $R$ and $C$ are strongly connected, which is commonly assumed in existing literature.

Then, the locally differentially private distributed online stochastic optimization algorithm for \eqref{problem2} is summarized in Algorithm \ref{algorithm1}. 
\begin{algorithm}[htp]
  \caption{Quantization-Based Locally Differentially Private Distributed Online Stochastic Optimization}
  \label{algorithm1}
  \textbf{Parameters: }Iteration stepsize $\lambda_t = \frac{\lambda_0}{(t+1)^\nu}$, dynamic quantization stepsize $d_t^i = \frac{d_0^i}{(t+1)^{\varsigma^i}}$, with $\lambda_0, d_0^i > 0$, $\nu, \varsigma^i \in (\frac{1}{2}, 1)$, and $\max\limits_{i \in [m]}\varsigma^i < \nu$.

  \textbf{Initialization: }Each agent $i$ maintains three variables $\theta_t^i$, $\psi_t^i$, and $z_t^i$, which are initialized with random values $\theta_0^i, \psi_0^i \in \mathbb{R}^d$, and $z_0^i = e_i \in \mathbb{R}^m$, where $e_i$ denotes the $i$th standard basis vector.

  \textbf{Iterations: }For $t = 0, 1, \cdots, T-1$, each agent $i \in [m]$ does

  a) Acquire a new data-point $x_t^i \sim \Phi^i$, and compute the gradient $\nabla f_t^i(\theta_t^i) = \frac{1}{t+1}\sum_{k=0}^{t}\nabla l_i(\theta_t^i, x_k^i)$ by using all acquired data up to time $t$ and the current optimization variable $\theta_t^i$.

  b) Use the dynamic stochastic quantizer \eqref{quantizer} to compute $Q(\theta_t^i)$ and $Q(\psi_t^i)$, then push $C_{pi}Q(\psi_t^i)$ to agent $p \in \mathbb{N}_{C,i}^\text{out}$ and pull $Q(\theta_t^j)$ from agent $j \in \mathbb{N}_{R,i}^\text{in}$, where the subscript $R$ and $C$ in neighbor sets indicates the neighbors with respect to the graphs induced by these matrices.

  c) Update its states as follows
  \begin{align}
    \psi_{t+1}^i &= (1+C_{ii})\psi_t^i + \sum_{j \in \mathbb{N}_{C,i}^\text{in}}C_{ij}Q(\psi_t^j) + \lambda_t \nabla f_t^i(\theta_t^i), \nonumber\\
    \theta_{t+1}^i &= (1+R_{ii})\theta_t^i + \sum_{j \in \mathbb{N}_{R,i}^\text{in}}R_{ij}Q(\theta_t^j) - \frac{\psi_{t+1}^i - \psi_t^i}{m[z_t^i]_i}, \nonumber\\
    z_{t+1}^i &= z_t^i + \sum_{j \in \mathbb{N}_{R,i}^\text{in}}R_{ij}(z_t^j - z_t^i), \label{states_update}
  \end{align}
  where $[z_t^i]_i$ denotes the $i$th element of $z_t^i$.
\end{algorithm}

\setcounter{counter}{2}
\refstepcounter{counter}  
\label{lemma 3}
\noindent\textbf{Lemma 3 \cite{Pu_Push-pull}.}
Under Assumption \ref{assumption 4}, the matrix $I+R$ has a unique nonnegative left eigenvector $u^T$ (associated with eigenvalue 1) satisfying $u^T \mathbf{1} = m$, and the matrix $I+C$ has a unique nonnegative right eigenvector $v$ (associated with eigenvalue 1) satisfying $\mathbf{1}^T v = m$.

According to \cite[Lemma 3]{Pu_Push-pull}, the spectral radius of $\bar{R} \triangleq I + R - \frac{\mathbf{1}u^T}{m}$ equals $1 - |\nu_R| < 1$ where $\nu_R$ is an eigenvalue of $R$. A vector norm $\|\cdot\|_R$ can be defined such that the induced matrix norm $\|\bar{R}\|_R$ is arbitrarily close to the spectral radius of $\bar{R}$, and in particular $\|\bar{R}\|_R < 1$. Without loss of generality, we represent this norm as $\|\bar{R}\|_R = 1 - \rho_R < 1$.
Similarly, we know that the spectral radius of $\bar{C} \triangleq I + C - \frac{v\mathbf{1}^T}{m}$ is equal to $1- |\nu_C|<1$ where $\nu_C$ is an eigenvalue of $C$. We also define a vector norm $\|\cdot\|_C$ such that the induced matrix norm $\|\bar{C}\|_C$ is arbitrarily close to the special radius of $\bar{C}$, and in particular $\|\bar{C}\|_C < 1$. Without loss of generality, this norm is bounded as $\|\bar{C}\|_C = 1 - \rho_C < 1$.

\setcounter{counter}{3}
\refstepcounter{counter}  
\label{lemma 4}
\noindent\textbf{Lemma 4 \cite{Wang_Gradient}.}
Under Assumption \ref{assumption 4}, the variables $z_t^i$ in \eqref{states_update}, after scaled by $m$, converge to the left eigenvector $u^T$ of $I+R$ with a geometric rate, i.e., there exist $C_z>0$ and $P_z\in(0,1)$ such that $\big|\frac{1}{m[z_t^i]_i} -\frac{1}{u_i}\big|\leq C_zP_z^t$ holds for any $i\in[m]$ and $t\ge 0$, where $u_i$ is the $i$th element of $u$.

\subsection{Convergence Analysis}
Now, we proceed to analyze the convergence of Algorithm \ref{algorithm1}. Defining $\xi_t^i = Q(\psi_t^i) - \psi_t^i$, $\zeta_t^i = Q(\theta_t^i) - \theta_t^i$ and $g_t^i = \nabla f_t^i(\theta_t^i)$, the first two equations of \eqref{states_update} can be rewritten as
\begin{align*}
  \psi_{t+1}^i &= (1+C_{ii})\psi_t^i + \sum_{j \in \mathbb{N}_{C,i}^\text{in}}C_{ij}(\psi_t^j + \xi_t^j) + \lambda_t g_t^i, \nonumber \\
  \theta_{t+1}^i &= (1+R_{ii})\theta_t^i + \sum_{j \in \mathbb{N}_{R,i}^\text{in}}R_{ij}(\theta_t^j + \zeta_t^j) - \frac{\psi_{t+1}^i - \psi_t^i}{m[z_t^i]_i}. \nonumber
\end{align*}
Then, defining the following stacked vectors $\xi_{t,C} = ((\xi_{t,C}^1)^T, \cdots, (\xi_{t,C}^m)^T)^T$ with $\xi_{t,C}^i = \sum_{j \in \mathbb{N}_{C,i}^\text{in}}C_{ij}\xi_t^j$, $\zeta_{t,R} = ((\zeta_{t,R}^1)^T, \cdots, (\zeta_{t,R}^m)^T)^T$ with $\zeta_{t,R}^i = \sum_{j\in \mathbb{N}_{R,i}^\text{in}}R_{ij}\zeta_t^j$,
$z_{t,m} = \text{diag}\{m[z_t^1]_1, \cdots, m[z_t^m]_m\} \otimes I_d$, $\psi_t=((\psi_t^1)^T,\cdots$,
$(\psi_t^m)^T)^T$, $\theta_t=((\theta_t^1)^T,\cdots,(\theta_t^m)^T)^T$, $z_t=((z_t^1)^T,\cdots$,
$(z_t^m)^T)^T$, and $g_t=((g_t^1)^T,\cdots,(g_t^m)^T)^T$, \eqref{states_update} in Algorithm \ref{algorithm1} can be written in the following compact form
\begin{align}
  \psi_{t+1} &= \hat{C}\psi_t + \xi_{t,C} + \lambda_t g_t, \nonumber\\
  \theta_{t+1} &= \hat{R} \theta_t + \zeta_{t,R} - z_{t,m}^{-1}(\psi_{t+1} - \psi_t), \nonumber\\
  z_{t+1} &= \hat{R}z_t, \label{vector_form}
\end{align}
where $\hat{R} = (I + R) \otimes I_d$ and $\hat{C}=(I +C) \otimes I_d$. For the convergence analysis, the following average vectors are defined as $\bar{\theta}_t = \frac{(u^T \otimes I_d)\theta_t}{m}$, $\bar{\psi}_t = \frac{(\mathbf{1}^T \otimes I_d)\psi_t}{m}$, $\bar{g}_t = \frac{(\mathbf{1}^T \otimes I_d)g_t}{m}$, $\bar{\xi}_{t,C} = \frac{(\mathbf{1}^T \otimes I_d)\xi_{t,C}}{m}$, $\bar{\zeta}_{t,R} = \frac{(u^T \otimes I_d)\zeta_{t,R}}{m}$, where $u$ is given in Lemma \ref{lemma 3}. Then, it follows from \eqref{vector_form} that $\bar{\psi}_{t+1} = \frac{1}{m}(\mathbf{1}^T \otimes I_d)((I + C) \otimes I_d)\psi_t + \bar{\xi}_{t,c} + \lambda_t \bar{g}_t$, $\bar{\theta}_{t+1} = \frac{1}{m}(u^T \otimes I_d)((I+R) \otimes I_d)\theta_t + \bar{\zeta}_{t,R} - \frac{(u^T \otimes I_d)}{m}U^{-1}(\psi_{t+1}-\psi_t) - \frac{(u^T \otimes I_d)}{m}(z_{t,m}^{-1}-U^{-1})(\psi_{t+1}-\psi_t)$,
where $U=\text{diag}(u)\otimes I_d$. Using the relationships $(\mathbf{1}^T \otimes I_d)((I + C) \otimes I_d) = \mathbf{1}^T \otimes I_d$, $(u^T \otimes I_d)((I+R) \otimes I_d) =u^T \otimes I_d$, and $(u^T \otimes I_d)U^{-1} = \mathbf{1}^T \otimes I_d$, one gets that
\begin{align}
  \bar{\psi}_{t+1} &= \bar{\psi_t} + \bar{\xi}_{t,C} + \lambda_t \bar{g}_t, \nonumber\\
  \bar{\theta}_{t+1} &= \bar{\theta}_t + \bar{\zeta}_{t,R} - \frac{\mathbf{1}^T \otimes I_d}{m}(\psi_{t+1}-\psi_t) \nonumber\\
  &\quad - \frac{(u^T \otimes I_d)}{m}(z_{t,m}^{-1}-U^{-1})(\psi_{t+1}-\psi_t). \label{average_state_1}
\end{align}
From Definition \ref{definition 1}, the vector norms $\|\cdot\|_R$ of $\theta_t$ and $\|\cdot\|_C$ of $\psi_t$ can be computed as
\begin{align*}
  \|\theta_t\|_R &= \|(\|\theta_t(1)\|_R, \|\theta_t(2)\|_R, \cdots, \|\theta_t(d)\|_R)\|_2, \nonumber\\
  \|\psi_t\|_C &= \|(\|\psi_t(1)\|_C, \|\psi_t(2)\|_C, \cdots, \|\psi_t(d)\|_C)\|_2,
\end{align*}
where $\theta_t(\ell)$ and $\psi_t(\ell)$ denote the $\ell$th column of $\theta_t$ and $\psi_t$, respectively, for $1 \leq \ell \leq d$. In the following, $\|\theta_t - \boldsymbol{\bar{\theta}_t}\|_R$ (with $\boldsymbol{\bar{\theta}_t} = \boldsymbol{1}_m \otimes \bar{\theta}_t$) is utilized to measure the distance between all $\theta_t^i$ and their average $\bar{\theta}_t$, and $\|\psi_t - \text{diag}(\boldsymbol{v})\boldsymbol{\bar{\psi}_t}\|_C$ (with $\boldsymbol{v} =\boldsymbol{1}_d \otimes v$ and $\boldsymbol{\bar{\psi}_t} = \boldsymbol{1}_m \otimes \bar{\psi}_t$) is used to measure the distance between all $\psi_t^i$ and their $v$-weighted average $\bar{\psi}_t$, where $v$ is given in Lemma \ref{lemma 3}. Next, a convergence result is provided for Algorithm \ref{algorithm1}.

\setcounter{counter}{0}
\refstepcounter{counter}  
\label{theorem 1}
\noindent\textbf{Theorem 1.}
For Algorithm \ref{algorithm1}, if Assumptions \ref{assumption 1}-\ref{assumption 4} are satisfied, then the following statements hold for each agent $i\in[m]$:
\\
\indent\setlength{\parindent}{1em}
a) $\lim_{t \rightarrow \infty}F(\bar{\theta_t})$ exists a.s. and $\lim_{t \rightarrow \infty} \|\theta_t^i - \bar{\theta}_t\| = \lim_{t \rightarrow \infty} \|\psi_t^i - v_i\bar{\psi}_t\| = 0$ a.s., where $v_i$ is the $i$th element of $v$.
\\
\indent\setlength{\parindent}{1em}
b) $\liminf_{t \rightarrow \infty}\|\nabla F(\bar{\theta}_t)\|=0$ holds a.s.. Moreover, every accumulation point of $\{\bar{\theta}_t\}$ is an optimal solution a.s., and $\lim_{t \rightarrow \infty}F(\theta_t^i)=F^*$ a.s..

\noindent\textbf{Proof.}
  See Appendix.

\noindent\textbf{Remark 2.}
Noting that the evolution of $\theta_t^i$ to the optimal solution satisfies the conditions in Lemmas \ref{lemma 1} and \ref{lemma 2}, we can utilize these lemmas to analyze the convergence speed. Specifically, it follows from \eqref{lim psi} and \eqref{lim theta} that $\|\psi_t^i-v_i\bar{\psi}_t\|^2$ and $\|\theta_t^i-\bar{\theta}_t\|^2$ decay to zero with a rate no slower than $\mathcal{O}(\frac{1}{t})$. Moreover, one can deduce from \eqref{nabla F} that $\|\nabla F(\bar{\theta}_t)\|^2$ decays to zero with a rate no slower than $\mathcal{O}(\frac{1}{t^{1-\nu}})$.
Obviously, the above convergence rates are independent of the quantization stepsize $d_t$, which implies that our algorithm with the dynamic stochastic quantizer does not pay a price in convergence rates.

\subsection{Privacy Analysis}
In this subsection, we prove that Algorithm \ref{algorithm1} can ensure $(0,\delta^i)$-LDP for all $i\in[m]$, even over an infinite number of iterations.
For a distributed online stochastic optimization algorithm, let $\mathcal{A}_i(\mathcal{D}_t^i,\theta_t^{-i})$ denote the implementation (i.e., internal states) of the algorithm by agent $i$, with $\mathcal{D}_t^i$ being the local dataset at iteration $t$ and $\theta_t^{-i}$ being all received messages. In Algorithm $1$, this implementation takes the specific form $\mathcal{A}_i(\mathcal{D}_{t}^i,\theta_{t}^{-i})=((\theta_{t+1}^i)^T,(\psi_{t+1}^i)^T,(z_{t+1}^i)^T )^T$.
Following \cite{Chen_Locally}, the sensitivity of the algorithm is defined as follows.

\setcounter{counter}{3}
\refstepcounter{counter}  
\label{definition 4}
\noindent\textbf{Definition 4 (Sensitivity).}
For any two adjacent datasets $\mathcal{D}_t^i$ and $\mathcal{D}_t'^{i}$ of agent $i\in [m]$ at each iteration $t$, the sensitivity of Algorithm \ref{algorithm1} is defined as
\begin{equation}
  \Delta_{t+1}^i=\max_{\text{Adj}(\mathcal{D}_t^i,\mathcal{D}_t'^{i})}\|\mathcal{A}_i(\mathcal{D}_t^i,\theta_t^{-i})-\mathcal{A}_i(\mathcal{D}_t'^{i},\theta_t^{-i})\|_1. \label{definition_sensitivity}
\end{equation}
\indent\setlength{\parindent}{1em}
According to Definition \ref{definition 4}, the sensitivity of Algorithm \ref{algorithm1} can be divided into three elements, i.e., $\Delta_{t,\theta}^i$, $\Delta_{t,\psi}^i$ and $\Delta_{t,z}^i$, which correspond to the three shared messages $\theta_t^i$, $\psi_t^i$ and $z_t^i$, respectively.
Noting that $z_t^i$ contains no sensitive information, it can be easily obtained that $\Delta_{t,z}^i = 0$.
For privacy analysis, we state the following lemmas.

\setcounter{counter}{4}
\refstepcounter{counter}  
\label{lemma 5}
\noindent\textbf{Lemma 5.}
Denote $\{v_t\}$ as a nonnegative sequence. If there exists a sequence $\eta_t=\frac{\eta_0}{(t+1)^n}$ with $\eta_0>0$, $n>0$ such that $v_{t+1}\leq (1-\varkappa)v_t+\eta_t$ holds for all $\varkappa\in(0,1)$ and $t\ge T_1$, where $T_1\in\mathbb{N}$. Then we obtain $v_t\leq \Xi\eta_t$ for all $t\ge T_1$, where $\Xi=\Big(\frac{4n}{e\text{ln}(\frac{2}{2-\varkappa})}\Big)^n \Big(\frac{v_{T_1}(1-\varkappa)^{1-T_1}}{\eta_0} +\frac{2}{\varkappa}\Big)$.

\noindent\textbf{Proof.}
  The lemma can be obtained following the same line of \cite[Lemma 4]{Chen_Local}.

\setcounter{counter}{5}
\refstepcounter{counter}  
\label{lemma 6}
\noindent\textbf{Lemma 6.}
Consider the dynamic stochastic quantizer \eqref{quantizer} with input $y\in\mathbb{R}$ and output $q=\mathcal{Q}(y)$. If $|y-y'|<d_t$, for any Borel-measurable set $\tau\subseteq\text{Range}(q)$, we have
\begin{equation}
  |\mathbb{P}(q\in\tau|y) -\mathbb{P}(q'\in\tau|y')|\leq \frac{|y-y'|}{d_t}.
\end{equation}
\noindent\textbf{Proof.}
  Since $|y-y'|<d_t$, without loss of generality, we consider three cases as
  \\
  \indent\setlength{\parindent}{1em}
  Case 1: $y, y'\in(nd_t,(n+1)d_t]$;
  \\
  \indent\setlength{\parindent}{1em}
  Case 2: $y\in(nd_t,(n+1)d_t]$, $y'\in((n+1)d_t,(n+2)d_t]$;
  \\
  \indent\setlength{\parindent}{1em}
  Case 3: $y\in(nd_t,(n+1)d_t]$, $y'\in((n-1)d_t,nd_t]$.
  \\
  \indent\setlength{\parindent}{1em}
  (Case 1) Define $\tau = \tau_1^* \cup \tau_1$ where $\tau_1^*=\{\emptyset,\{nd_t\},\{(n+1)d_t\},\{nd_t,(n+1)d_t\}\}$ and $\tau_1=\{N_1\cdot d_t|N_1\subseteq\mathbb{Z},\ N_1\cap\{n,n+1\}=\emptyset\}$. Since $\mathbb{P}(q\in\tau_1|y)=\mathbb{P}(q'\in\tau_1|y')=0$, we have $\mathbb{P}(q\in\tau|y)=\mathbb{P}(q\in\tau_1^*|y)$, $\mathbb{P}(q'\in\tau|y')=\mathbb{P}(q'\in\tau_1^*|y')$, i.e., $|\mathbb{P}(q\in\tau|y) -\mathbb{P}(q'\in\tau|y')| = |\mathbb{P}(q\in\tau_1^*|y) -\mathbb{P}(q'\in\tau_1^*|y')|$.
  \\
  \indent\setlength{\parindent}{1em}
  Firstly, for $\tau=\emptyset$ and $\tau=\{nd_t,(n+1)d_t\}$, it can be derived that
  \begin{equation}
    |\mathbb{P}(q\in\tau|y) -\mathbb{P}(q'\in\tau|y')|=0.
  \end{equation}
  Next, it follows from $\tau=\{nd_t\}$ that
  \begin{equation*}
    \begin{aligned}
      &|\mathbb{P}(q=nd_t|y) -\mathbb{P}(q'=nd_t|y')| \nonumber\\
      &= \left|\left(1-\frac{y-nd_t}{d_t}\right)-\left(1-\frac{y'-nd_t}{d_t}\right)\right| = \frac{|y-y'|}{d_t}.
    \end{aligned}
  \end{equation*}
  Finally, for $\tau=\{(n+1)d_t\}$, we obtain
  \begin{align}
    &|\mathbb{P}(q=(n+1)d_t|y) -\mathbb{P}(q'=(n+1)d_t|y')| \nonumber\\
    &= \left|\frac{y-nd_t}{d_t}-\frac{y'-nd_t}{d_t}\right| = \frac{|y-y'|}{d_t}.
  \end{align}
  \\
  \indent\setlength{\parindent}{1em}
  (Case 2) Define $\tau=\tau_2^*\cup\tau_2$ where $\tau_2^*=\{\emptyset,\{nd_t\},\{(n+1)d_t\},\{(n+2)d_t\},\{nd_t,(n+1)d_t\},\{nd_t,(n+2)d_t\},\{(n+1)d_t,(n+2)d_t\},\{nd_t,(n+1)d_t,(n+2)d_t\}\}$ and $\tau_2=\{N_2\cdot d_t|N_2\subseteq\mathbb{Z},N_2\cap\{n,n+1,n+2\}=\emptyset\}$. According to $\mathbb{P}(q\in\tau_2|y)=\mathbb{P}(q'\in\tau_2|y')=0$, we have $\mathbb{P}(q\in\tau|y)=\mathbb{P}(q\in\tau_2^*|y)$, $\mathbb{P}(q'\in\tau|y')=\mathbb{P}(q'\in\tau_2^*|y')$, i.e., $|\mathbb{P}(q\in\tau|y) -\mathbb{P}(q'\in\tau|y')| = |\mathbb{P}(q\in\tau_2^*|y) -\mathbb{P}(q'\in\tau_2^*|y')|$.
  \\
  \indent\setlength{\parindent}{1em}
  Firstly, for $\tau=\emptyset$ and $\tau=\{nd_t,(n+1)d_t,(n+2)d_t\}$, it is straightforward that
  \begin{equation}
    |\mathbb{P}(q\in\tau|y) -\mathbb{P}(q'\in\tau|y')|=0.
  \end{equation}
  Secondly, for $\tau=\{nd_t\}$, we arrive at
  \begin{align}
    &|\mathbb{P}(q=nd_t|y) -\mathbb{P}(q'=nd_t|y')| \nonumber\\
    &= \left|(1-\frac{y-nd_t}{d_t})-0\right| \leq \frac{|y-y'|}{d_t}.
  \end{align}
  The case of $\tau=\{(n+2)d_t\}$ can be obtained similarly.
  Thirdly, it follows from $\tau=\{(n+1)d_t\}$ that
  \begin{align}
    &|\mathbb{P}(q=(n+1)d_t|y) -\mathbb{P}(q'=(n+1)d_t|y')| \nonumber\\
    &= \left|\frac{y-nd_t}{d_t}-\left(1-\frac{y'-(n+1)d_t}{d_t}\right)\right| \nonumber\\
    &\leq \frac{(n+1)d_t-y+y'-(n+1)d_t}{d_t} = \frac{|y-y'|}{d_t}.
  \end{align}
  Finally, for $\tau=\{nd_t,(n+1)d_t\}$, we obtain
  \begin{equation*}
    \begin{aligned}
      &|\mathbb{P}(q\in\{nd_t,(n+1)d_t\}|y) -\mathbb{P}(q'\in\{nd_t,(n+1)d_t\}|y')| \nonumber\\
      &=|(1-\mathbb{P}(q=(n+2)d_t|y))-(1-\mathbb{P}(q'=(n+2)d_t|y'))| \nonumber\\
      &=|\mathbb{P}(q=(n+2)d_t|y)-\mathbb{P}(q'=(n+2)d_t|y')| \leq \frac{|y-y'|}{d_t}.
    \end{aligned}
  \end{equation*}
  The result can also be obtained for the cases where $\tau=\{nd_t,(n+2)d_t\}$ and $\tau=\{(n+1)d_t,(n+2)d_t\}$, and Case 3 yields the same conclusions as Case 2. The proof of Lemma 6 is complete.

With the help of Lemmas \ref{lemma 5} and \ref{lemma 6}, the following results can be obtained.

\setcounter{counter}{1}
\refstepcounter{counter}  
\label{theorem 2}
\noindent\textbf{Theorem 2.}
For Algorithm \ref{algorithm1}, if Assumptions \ref{assumption 1}-\ref{assumption 4} are satisfied, then the following statements hold  for each agent $i\in[m]$:
\\
\indent\setlength{\parindent}{1em}
1) For any finite number of iterations $T$, Algorithm \ref{algorithm1} is $(0,\delta^i)$-LDP if the initial quantization step $d_0^i$ satisfies $d_0^i\ge\sum_{t=0}^{T-1}(\varrho_{t,\psi}^i +\varrho_{t,\theta}^i)(t+1)^{\varsigma^i}$ with $\varrho_{t,\psi}^i =2d_l\sum_{p=0}^{t-1}(1-|C_{ii}|)^p\lambda_{t-1-p}$ and $\varrho_{t,\theta}^i =\sum_{p=0}^{t-1}(1-|R_{ii}|)^{t-p-1}(C_zP_z^{p} +\frac{1}{|u_i|})(\varrho_{p+1,\psi}^i +\varrho_{p,\psi}^i)$.
\\
\indent\setlength{\parindent}{1em}
2) Algorithm \ref{algorithm1} can ensure $(0,\delta^i)$-LDP over infinite iterations if the initial quantization step $d_0^i$ satisfies $d_0^i\ge\sum_{t=0}^{T-1}\frac{C_3(C_z+C_0)}{(t+1)^{1+\nu-\varsigma^i}}$, where $C_0$ and $C_3$ are given in \eqref{Delta_theta_psi} and \eqref{Delta_t^i_2}, respectively.

\noindent\textbf{Proof.}
  1) For any finite number of iterations $T$, we first compute the sensitivity $\Delta_t^i$ for Algorithm \ref{algorithm1}.
  According to the definition of sensitivity in \eqref{definition_sensitivity}, we have $\mathcal{Q}(\theta_t^j)=\mathcal{Q}(\theta_t'^{j})$ and $\mathcal{Q}(\psi_t^j)=\mathcal{Q}(\psi_t'^j)$ for any $t\ge 0$ and $j\in\mathbb{N}_i^\text{in}$. Since the adjacent datasets $\mathcal{D}_t^i=\{x_1^i,\cdots,x_k^i,\cdots,x_t^i\}$ and $\mathcal{D}_t'^i=\{x_1'^i,\cdots,x_k'^i,\cdots,x_t'^i\}$ are different only at time instant $k$, we have $x_p^i=x_p'^i$ for all $p\ne k$. Simultaneously, due to the difference between $D_t^i$ and $D_t'^i$, the loss function will change at time $k$, i.e., $l_i(\theta,x_k^i)\ne l_i(\theta,x_k'^i)$, and we have $\theta_t^i\ne \theta_t'^i$, $\psi_t^i\ne \psi_t'^i$, for $t>k$. Hence, $
    \|\psi_{t+1}^i -\psi_{t+1}'^i\|_1 = \|(1+C_{ii})(\psi_t^i -\psi_t'^i) + \frac{\lambda_t}{t+1}\sum_{p=0}^{t}(\nabla l_i(\theta_t^i, x_p^i)-\nabla l_i(\theta_t'^i, x_p'^i)\|_1$, and the sensitivity $\Delta_{t+1,\psi}^i$ satisfies
  \begin{align}
    \Delta_{t+1,\psi}^i &\leq \frac{\lambda_t}{t+1}\sum_{p=0}^{t}\|\nabla l_i(\theta_t^i, x_p^i)-\nabla l_i(\theta_t'^i, x_p'^i)\|_1\nonumber \\
    &\quad + (1-|C_{ii}|)\Delta_{t,\psi}^i.\label{psi-psi'}
  \end{align}
  By using Assumption \ref{assumption 3} and the property $\Delta_{0,\psi}^i=0$, we iterate the above inequality from $0$ to $t-1$ and obtain
  \begin{equation}
    \Delta_{t,\psi}^i\leq 2d_l\sum_{p=0}^{t-1}(1-|C_{ii}|)^p\lambda_{t-1-p}. \label{Delta_psi_varrho}
  \end{equation}
  Similarly, the following relation holds $\|\theta_{t+1}^i -\theta_{t+1}'^i\|_1 = \|(1+R_{ii})(\theta_t^i -\theta_t'^i) 
  -\frac{(\psi_{t+1}^i -\psi_{t+1}'^i)}{m[z_t^i]_i} +\frac{(\psi_t^i -\psi_t'^i)}{m[z_t^i]_i}\|_1$.
  By using Lemma \ref{lemma 4}, the sensitivity $\Delta_{t+1,\theta}^i$ satisfies $\Delta_{t+1,\theta}^i \leq  (C_zP_z^t +\frac{1}{|u_i|})(\Delta_{t+1,\psi}^i +\Delta_{t,\psi}^i) +(1-|R_{ii}|)\Delta_{t,\theta}$.
  Under the property $\Delta_{0,\theta}^i=0$, we iterate the above inequality from $0$ to $t-1$ and obtain $\Delta_{t,\theta}^i \leq \sum_{p=0}^{t-1}(1-|R_{ii}|)^{t-p-1}(C_zP_z^{p} +\frac{1}{|u_i|})(\Delta_{p+1,\psi}^i +\Delta_{p,\psi}^i)$. Since $\mathcal{A}_i(\mathcal{D}_{t-1}^i,\theta_{t-1}^{-i})=((\theta_t^i)^T$, $(\psi_t^i)^T,(z_t^i)^T )^T$ and $\Delta_{t,z}^i=0$, it can be obtained that
  \begin{equation}
    \Delta_t^i = \Delta_{t,\theta}^i +\Delta_{t,\psi}^i \leq \varrho_{t,\psi}^i + \varrho_{t,\theta}^i, \label{Delta_t^i_1}
  \end{equation}
  where $\varrho_{t,\psi}^i$ and $\varrho_{t,\theta}^i$ are given in the theorem statement.
  \\
  \indent\setlength{\parindent}{1em}
  From \eqref{Delta_t^i_1} and the condition on $d_0^i$ in Theorem \ref{theorem 2}, it follows that $\sum_{t=0}^{T-1}\frac{\Delta_t^i}{d_t^i} \leq 1$ for any finite number of iterations $T$, and that $\Delta_t^i < d_t^i$ holds for all $t\ge 0$, $i\in[m]$.
  \\
  \indent\setlength{\parindent}{1em}
  By Definition \ref{definition 3}, we have $\mathcal{M}_t(\mathcal{D}^i) =\mathcal{Q}(\mathcal{A}_i(\mathcal{D}_{t-1}^i,\theta_{t-1}^{-i}))$ $=q_t$ and $\mathcal{M}_t(\mathcal{D}'^i) =\mathcal{Q}(\mathcal{A}_i(\mathcal{D}_{t-1}'^i,\theta_{t-1}^{-i})) =q_t'$, where
  \begin{align*}
    &q_t=(\mathcal{Q}([\theta_t^i]_1),\cdots,\mathcal{Q}([\theta_t^i]_d),\mathcal{Q}([\psi_t^i]_1),\cdots,\mathcal{Q}([\psi_t^i]_d))^T, \\
    &q_t'=(\mathcal{Q}([\theta_t'^i]_1),\cdots,\mathcal{Q}([\theta_t'^i]_d),\mathcal{Q}([\psi_t'^i]_1),\cdots,\mathcal{Q}([\psi_t'^i]_d))^T.
  \end{align*}
  Note that $[q_t]_{\ell=1,\cdots,2d}$ denotes the $i$th element of $q_t$, which can be seen as a single LDP mechanism $\mathcal{M}_t^\ell(\mathcal{D}^i)$. Given $\Delta_t^i < d_t^i$, it follows from Lemma \ref{lemma 6} that
  \begin{align}
    &|\mathbb{P}[\mathcal{M}_t^\ell(\mathcal{D}^i)\in\mathcal{S}_t^\ell] -\mathbb{P}[\mathcal{M}_t^\ell(\mathcal{D}'^i)\in\mathcal{S}_t^\ell]| \nonumber\\
    &\quad= \Big|\mathbb{P}\left[[q_t]_\ell\in\mathcal{S}_t^\ell \big|[\mathcal{A}_i(\mathcal{D}_{t-1}^i,\theta_{t-1}^{-i})]_\ell\right] \nonumber\\
    &\qquad -\mathbb{P}\left[[q_t']_\ell\in\mathcal{S}_t^\ell \big|[\mathcal{A}_i(\mathcal{D}_{t-1}'^i,\theta_{t-1}^{-i})]_\ell\right]\Big| \nonumber\\
    &\quad\leq \frac{|[\mathcal{A}_i(\mathcal{D}_{t-1}^i,\theta_{t-1}^{-i})]_\ell -[\mathcal{A}_i(\mathcal{D}_{t-1}'^i,\theta_{t-1}^{-i})]_\ell|}{d_t^i}, \label{P_ell}
  \end{align}
  where $\mathcal{S}_t^\ell \subseteq \text{Range}(\mathcal{M}_t^\ell)$.
  Defining $|[\mathcal{A}_i(\mathcal{D}_{t-1}^i,\theta_{t-1}^{-i})]_\ell -[\mathcal{A}_i(\mathcal{D}_{t-1}'^i,\theta_{t-1}^{-i})]_\ell|/d_t^i=\delta_{t,\ell}^i$, it follows from \eqref{P_ell} that the mechanism $\mathcal{M}_t^\ell(\mathcal{D}^i)$ can ensure $(0,\delta_{t,\ell}^i)$-LDP for any $\ell=1,\cdots,2d$, where $\delta_{t,\ell}^i\in[0,1)$.
  Then, by using the composition theorem for DP (which directly applies to LDP) in \cite{Dwork_The}, we obtain that the composition mechanism $\mathcal{M}_t(\mathcal{D}^i)$ can guarantee $(0,\sum_{\ell=1}^{2d}\delta_{t,\ell}^i)$-LDP, i.e.,
  \begin{equation}
    |\mathbb{P}[\mathcal{M}_t(\mathcal{D}^i)\in\mathcal{S}_t] -\mathbb{P}[\mathcal{M}_t(\mathcal{D}'^i)\in\mathcal{S}_t]| \leq \sum_{\ell=1}^{2d}\delta_{t,\ell}^i
  \end{equation}
  for any Borel-measurable set $\mathcal{S}_t \subseteq \text{Range}(\mathcal{M}_t)$. By defining $\delta_t^i=\frac{\Delta_t^i}{d_t^i}$, since $\sum_{\ell=1}^{2d}\delta_{t,\ell}^i \leq\frac{\Delta_t^i}{d_t^i}$ and $\delta_t^i < 1$, $\mathcal{M}_t(\mathcal{D}^i)$ achieves $\left(0,\delta_t^i\right)$-LDP for any iteration $t=0,\cdots,T-1$.
  \\
  \indent\setlength{\parindent}{1em}
  Since $\mathcal{M}(\mathcal{D}^i)=(\mathcal{M}_0(\mathcal{D}^i),\cdots,\mathcal{M}_{T-1}(\mathcal{D}^i))$, we use the composition theorem again, and obtain that the composition mechanism $\mathcal{M}(\mathcal{D}^i)$ is $(0,\sum_{t=0}^{T-1}\delta_t^i)$-LDP, i.e.,
  \begin{equation}
    |\mathbb{P}[\mathcal{M}(\mathcal{D}^i)\in\mathcal{S}] -\mathbb{P}[\mathcal{M}(\mathcal{D}'^i)\in\mathcal{S}]| \leq \sum_{t=0}^{T-1}\delta_t^i
  \end{equation}
  for any Borel-measurable set $\mathcal{S}\subseteq \text{Range}(\mathcal{M})$.
  Similarly, for any finite $T>0$, it follows from $\sum_{t=0}^{T-1}\delta_t^i \leq 1$ that the mechanism $\mathcal{M}(\mathcal{D}^i)$ is $(0,\delta^i)$-LDP by defining $\delta^i=\sum_{t=0}^{T-1}\delta_t^i$, i.e., Algorithm \ref{algorithm1} can ensure $(0,\delta^i)$-LDP for each agent $i$ for any finite number of iterations $T$.
 
  2) As $T$ tends to infinity, according to Assumptions \ref{assumption 2} and \ref{assumption 3}, \eqref{psi-psi'} can be rewritten as
  \begin{equation}
    \Delta_{t+1,\psi}^i \leq (1-|C_{ii}|)\Delta_{t,\psi}^i +\frac{\sqrt{n}Lt\lambda_t\Delta_{t,\theta}^i}{t+1} +\frac{2d_l\lambda_t}{t+1}.
    \label{Delta_psi_2}
  \end{equation}
  Substituting \eqref{Delta_psi_2} into the upper bound of $\Delta_{t+1,\theta}^i$ yields
  \begin{align}
    \Delta_{t+1,\theta}^i &\leq \left(1-|R_{ii}|+\frac{\sqrt{n}LC_zP_z^t\lambda_t}{t+1}\right)\Delta_{t,\theta}^i \nonumber\\
    &\quad +(2-|C_{ii}|)C_zP_z^t\Delta_{t,\psi}^i+\frac{1}{|u_i|}\Delta_{t,\psi}^i\nonumber\\
    &\quad +\frac{2d_lC_zP_z^t\lambda_t}{t+1} +\frac{1}{|u_i|}\Delta_{t+1,\psi}^i . \label{Delta_theta_2}
  \end{align}
  By choosing positive constants $C_1<\min\left\{\frac{|R_{ii}|}{2},\frac{|C_{ii}|}{2}\right\}$ and $C_0\ge \frac{4-2C_1}{|u_i|(|C_{ii}|-2C_1)}$, it can be deduced from \eqref{Delta_psi_2} and \eqref{Delta_theta_2} that
  \begin{align}
    &\Delta_{t+1,\theta}^i +\left(C_0 -\frac{1}{|u_i|}\right)\Delta_{t+1,\psi}^i \nonumber\\
    &\leq \left(1-|R_{ii}| +\frac{\sqrt{n}LC_z(tP_z^t\lambda_t)}{t+1} +\frac{C_0\sqrt{n}Lt\lambda_t}{t+1}\right)\Delta_{t,\theta}^i \nonumber\\
    &\quad +\left((2-|C_{ii}|)C_zP_z^t +C_0(1-|C_{ii}|) +\frac{1}{|u_i|}\right)\Delta_{t,\psi}^i \nonumber\\
    &\quad +\frac{2d_lC_zP_z^t\lambda_t +2d_lC_0\lambda_t}{t+1}. \label{Delta_theta_psi}
  \end{align}
  Since $\lambda_t$ and $P_z^t$ are both decaying, there must exists some $T_0>0$ such that $\frac{|R_{ii}|}{2}\ge \frac{\sqrt{n}LC_z(tP_z^t\lambda_t)}{t+1} +\frac{C_0\sqrt{n}Lt\lambda_t}{t+1}$ and $\frac{C_0|C_{ii}|}{2}\ge (2-|C_{ii}|)C_zP_z^t$ hold for all $t \ge T_0$. Hence, we always have $\Delta_{t+1,\theta}^i +(C_0 -\frac{1}{|u_i|})\Delta_{t+1,\psi}^i \leq (1-C_1)(\Delta_{t,\theta}^i +(C_0 -\frac{1}{|u_i|})\Delta_{t,\psi}^i) +\frac{2d_l\lambda_0(C_z +C_0)}{(t+1)^{1+\nu}}$ for $t\ge T_0$.
  Then, it follows from Lemma \ref{lemma 5} that for all $t\ge 0$,
  \begin{equation}
    \Delta_{t,\theta}^i +\left(C_0 -\frac{1}{|u_i|}\right)\Delta_{t,\psi}^i \leq C_2\frac{2d_l\lambda_0(C_z +C_0)}{(t+1)^{1+\nu}},
    \label{Delta_theta_psi_3}
  \end{equation}
  where the constant $C_2$ satisfies $C_2=\max_{i\in[m]} \Big\{(\frac{2}{C_1}+$
  $\frac{\varrho_{T_0}^i(1-C_1)^{1-T_0}}{2d_l\lambda_0(C_z +C_0)})(\frac{4(1+\nu)}{e\text{ln}(\frac{2}{2-C_1})})^{1+\nu}, \frac{\varrho_{T_0}^i(T_0+1)^{1+\nu}}{2d_l\lambda_0(C_z+C_0)} \Big\}$ with $\varrho_{T_0}^i = \varrho_{T_0,\theta}^i +\big(C_0 -\frac{1}{|u_i|}\big)\varrho_{T_0,\psi}^i$.
  \\
  \indent\setlength{\parindent}{1em}
  From \eqref{Delta_theta_psi_3}, for all $t>0$, the following relations hold $\Delta_{t,\theta}^i \leq C_2\frac{2d_l\lambda_0(C_z +C_0)}{(t+1)^{1+\nu}}$ and $\Delta_{t,\psi}^i \leq C_2\frac{2d_l\lambda_0(C_z +C_0)}{(C_0-\frac{1}{|u_i|})(t+1)^{1+\nu}}$.
  By defining $C_3=\frac{2C_2d_l\lambda_0((C_0+1)|u_i|-1)}{C_0|u_i|-1}$, it follows that
  \begin{equation}
    \Delta_t^i \leq C_3\frac{C_z+C_0}{(t+1)^{1+\nu}}.
    \label{Delta_t^i_2}
  \end{equation}
  \indent\setlength{\parindent}{1em}
  The rest of the proof is similar to the first part of Theorem \ref{theorem 2}, from which it follows that Algorithm \ref{algorithm1} remains $(0,\delta^i)$-LDP for each agent $i$ over infinite iterations.
  Furthermore, since $\delta^i=\sum_{t=0}^{T-1}\frac{\Delta_t^i}{d_t^i}$, a larger quantization stepsize $d_t^i$ leads to a smaller $\delta^i$, which implies stronger privacy protection. The proof of Theorem 2 is complete.

\noindent\textbf{Remark 3.}
Theorem \ref{theorem 2} establishes that Algorithm 1 is capable of ensuring LDP over an infinite number of iterations, which is different from \cite{Wang_Quantization} and \cite{Lin_Quantization} that only guarantee DP on a per-iteration basis.
The key idea is to collaboratively design the quantization stepsize and the iterative stepsize to ensure that $\delta_t^i$ is summable, thereby guaranteeing that the condition regarding $d_0^i$ in Theorem \ref{theorem 2} is achievable.

\noindent\textbf{Remark 4.}
In the degenerate case where the problem of locally differentially private distributed online stochastic optimization reduces to the standard DP framework, all agents share a unified privacy budget $\delta$ rather than individual $\delta^i$ for each agent $i$.
Consequently, Algorithm \ref{algorithm1} can achieve $(0, \delta)$-DP over an infinite number of iterations by appropriately adopting a uniform quantization stepsize $d_t$. This guarantee remains stronger than the per-iteration DP guarantees provided by \cite{Wang_Quantization} and \cite{Lin_Quantization}.

\section{Numerical Experiments}
\begin{figure}[htp]
  \centering
  \includegraphics[width=0.16\textwidth]{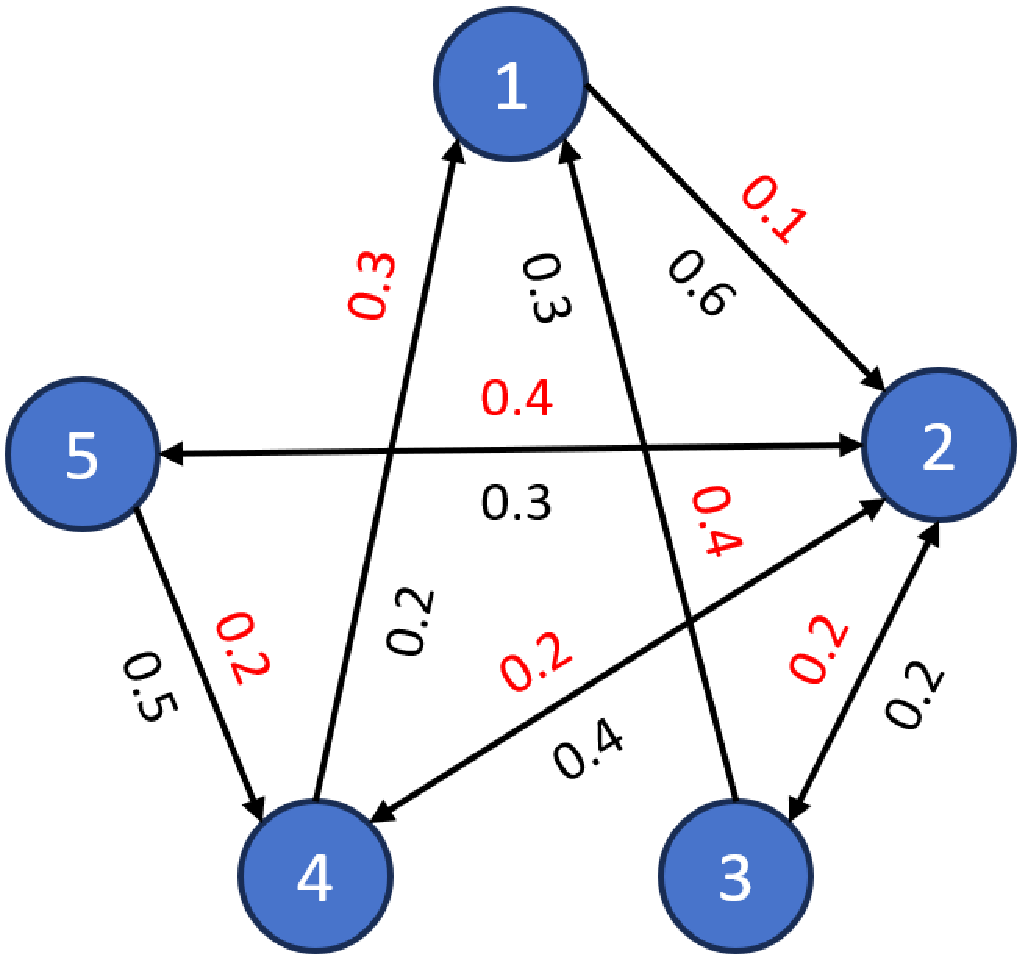}
  \caption{Topology structure of the directed graph}
  \label{figure_network}
\end{figure}

In this section, we verify the effectiveness of Algorithm \ref{algorithm1} through distributed online training of three classic machine learning experiments by employing the ``Mushrooms'', ``MNIST'', and ``BNCI2014001'' datasets. The communication topology is depicted in Fig. \ref{figure_network}, where the weights of matrices $R$ and $C$ are indicated in red and black, respectively.

\subsection{Logistic Regression on the ``Mushrooms” Dataset}
In the first experiment, we evaluate Algorithm \ref{algorithm1} on a logistic regression classification task by processing the ``Mushrooms'' dataset \cite{Dua_UCI} in an online manner, where data instances are acquired sequentially in real-time. In fact, each agent $i$ may have a distinct loss function $l_i$. In this case, the loss function in problem \eqref{problem2} is set as $l_i(\theta, x_k^i) = \frac{1}{N_k^i}\sum_{s=1}^{N_k^i} \text{log}\big(1 + e^{-b_t^i(s)a_t^i(s)^T\theta}\big)$ for all $i$, where $N_k^i$ is the number of sampling of agent $i$ at iteration $k$, while $a_t^i$ and $b_t^i$ denote the corresponding samples and labels, respectively.
\begin{figure}[htp]
  \centering
  \includegraphics[width=0.45\textwidth]{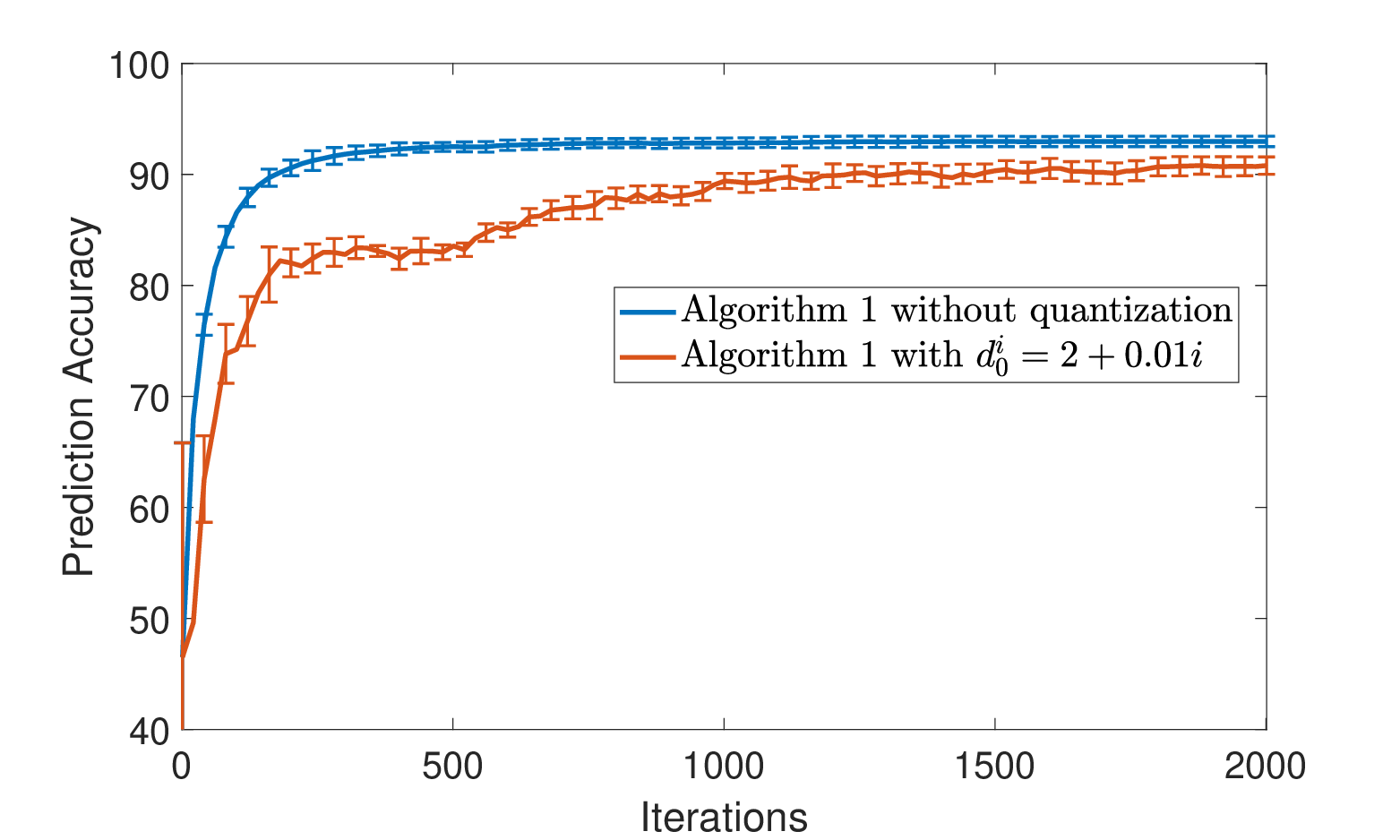}
  \caption{Comparison of the prediction accuracy of Algorithm 1 with and without quantization}
  \label{figure_mushrooms_acc}
\end{figure}
\begin{figure}[htp]
  \centering
  \includegraphics[width=0.45\textwidth]{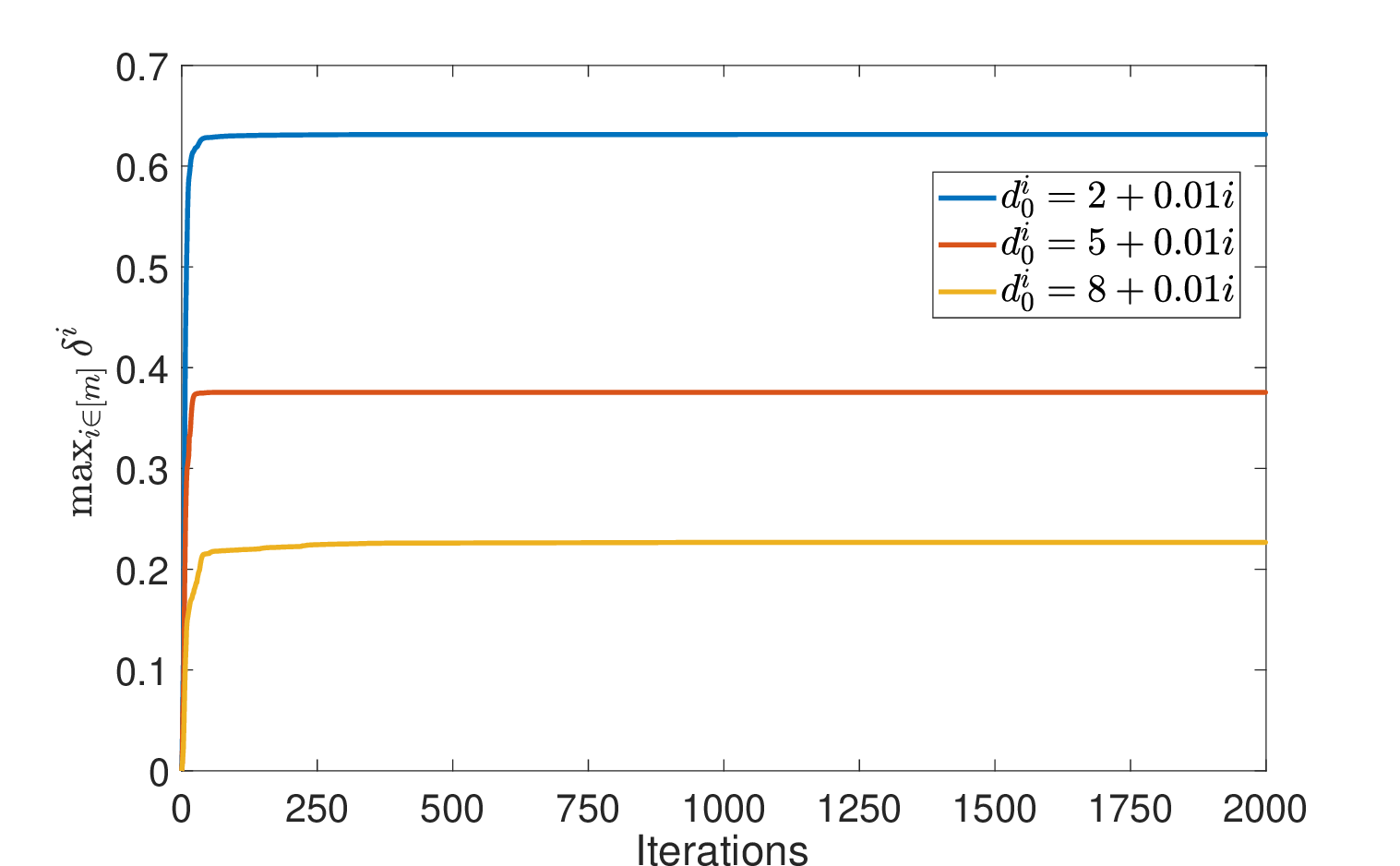}
  \caption{Comparison of the maximum $\delta^i$ under different quantization stepsizes}
  \label{figure_mushrooms_delta}
\end{figure}
\\
\indent\setlength{\parindent}{1em}
In each iteration, each agent randomly collects two labeled samples for training. The dynamic quantization stepsize and the iteration stepsize are configuerd as $d_t^i = \frac{d_0^i}{(t+1)^{\varsigma^i}}$ with $d_0^i = 2+0.01i$, $\varsigma^i = 0.6 + 0.01i$ and $\lambda_t = \frac{0.5}{(t+1)^{0.71}}$ for $i=1,2,\cdots,5$, respectively.
Prior to formal training, each agent sequentially receives 2000 random samples to serve as its test set. Once training begins, the prediction accuracy of the current logistic regression model is evaluated every 40 iterations using this test set.
The evolution of the average prediction accuracy and variance over 2000 runs are illustrated by the solid orange curve with error bars in Fig. \ref{figure_mushrooms_acc}. For comparison, the baseline is set as the non-quantized form of Algorithm \ref{algorithm1} which can be seen as a distributed online stochastic optimization algorithm without privacy protection.
From Fig. \ref{figure_mushrooms_acc}, Algorithm \ref{algorithm1} with quantization has a comparable convergence accuracy to that of the baseline.
This demonstrates the strong robustness of the proposed algorithm against errors induced by stochastic quantization, enabling it to achieve accurate convergence.
\\
\indent\setlength{\parindent}{1em}
To illustrate the privacy protection effect of Algorithm \ref{algorithm1}, we randomly modify a data point obtained from each agent and calculate the maximum value of $\delta^i$ accumulated across iterations. This metric serves as an indicator of the worst-case privacy level, and its variation curve is shown in Fig. \ref{figure_mushrooms_delta}.
As can be seen from the figure, $\max_{i\in[m]}\delta^i$ no longer changes rapidly after approximately 50 steps and stabilizes within the valid range, thereby ensuring a continuous privacy protection.
This is because the quantization stepsize and the iteration stepsize are jointly designed such that the decay rate of the quantization stepsize is slower than that of the sensitivity.
Furthermore, we also provide the evolution of $\max_{i\in[m]}\delta^i$ under $d_0^i=5+0.01i$ and $d_0^i=8+0.01i$ for $i=1,2,\cdots,5$. It can be observed that a larger quantization stepsize yields better privacy protection, which is consistent with our theoretical results.

\subsection{Image Recognition on the ``MNIST'' Dataset}
\begin{figure}[htp]
  \centering
  \includegraphics[width=0.48\textwidth]{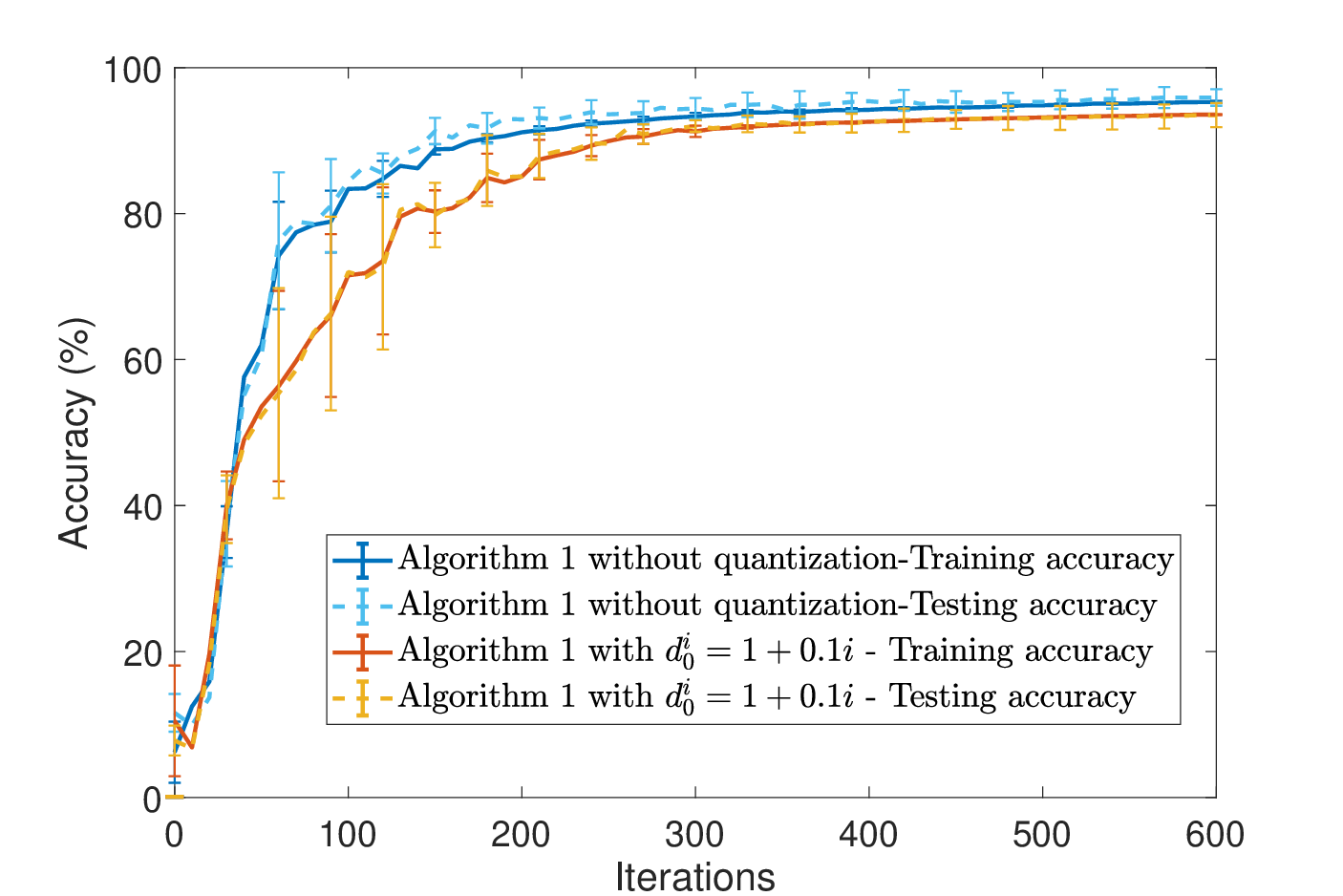}
  \caption{Comparison of the training and testing accuracy of Algorithm 1 with and without quantization}
  \label{figure_mnist_acc}
\end{figure}

In this experiment, we perform distributed online training of a convolutional neural network (CNN) on the “MNIST” dataset \cite{Lucun_Gradient-based}, a large benchmark
database of handwritten digits widely used for training and testing in machine learning. In the training, each agent has a local copy of the CNN. The CNN consists of four convolutional layers (the first two with 32 filters, the last two with 64 filters), each followed by a max pooling layer, and a final dense layer with 256 units.
Before the iteration begins, each agent stores 200 randomly selected images that arrive sequentially as the test set.
Once the formal iteration starts, each agent randomly selects 20 images per iteration for training.
After every 10 iterations, the training accuracy is computed using all the data collected so far, while the test accuracy is evaluated on the pre-stored test set.
To evaluate the proposed algorithm, the quantization stepsize and iteration stepsize are set as $d_t^i = \frac{1+0.1i}{(t+1)^{0.5+0.01i}}$ and $\lambda_t = \frac{0.1}{(t+1)^{0.61}}$ for $i=1,2,\cdots,5$, respectively. The evolution of the average training and testing accuracies over 600 iterations is illustrated by the solid blue and dashed light blue curves with error bars in Fig. \ref{figure_mnist_acc}. For comparison, we also set the non-private version of our algorithm as the baseline, whose accuracy curves are shown by the solid orange and dotted golden curves with error bars in Fig. \ref{figure_mnist_acc}. It is obvious that Algorithm \ref{algorithm1} with quantization achieves a convergence accuracy comparable to the baseline, even when training a complex CNN model.
\begin{figure}[htp]
  \centering
  \includegraphics[width=0.45\textwidth]{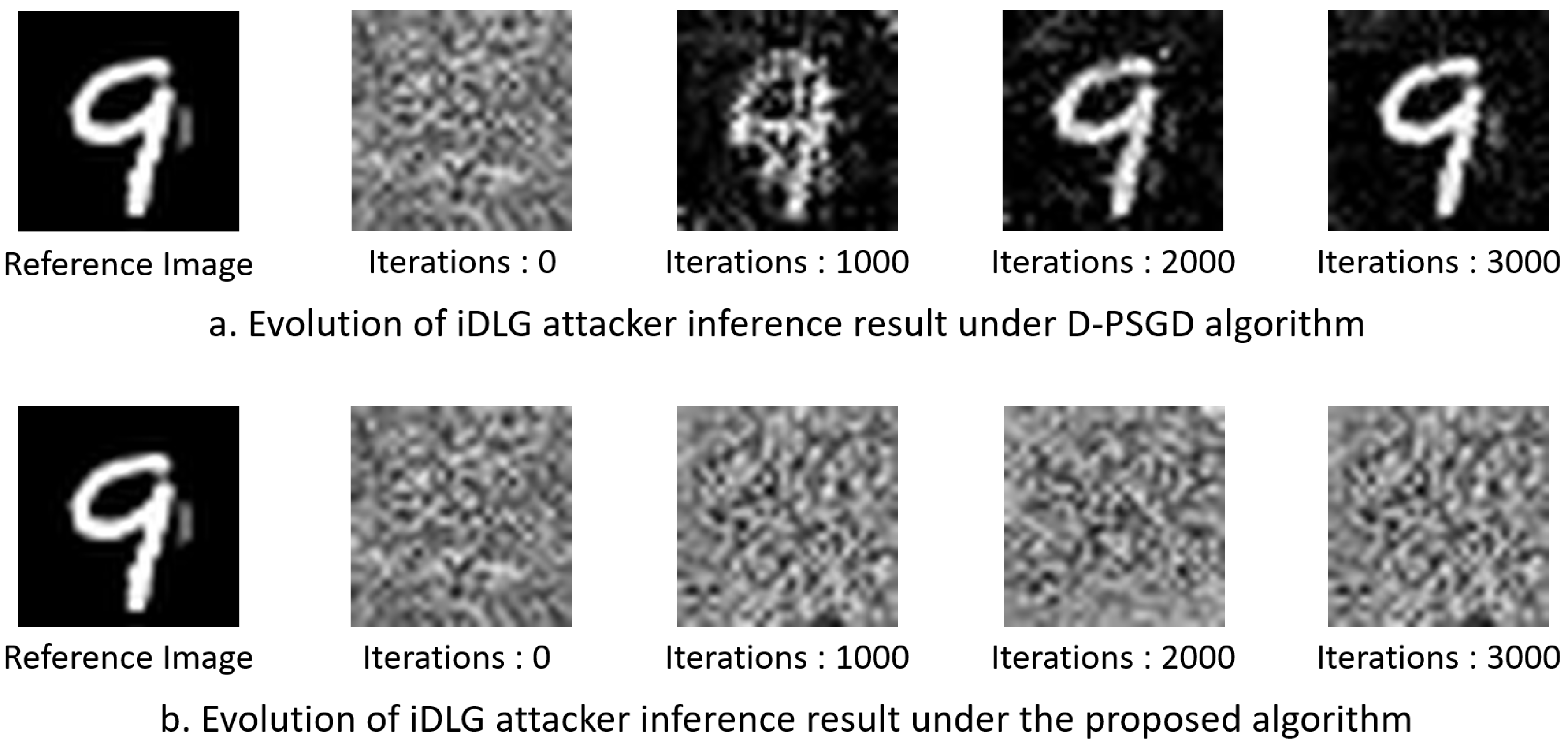}
  \caption{Comparison of iDLG attacker's image inference results under D-PSGD algorithm and our algorithm}
  \label{figure_idlg}
\end{figure}
\\
\indent\setlength{\parindent}{1em}
To verify the privacy protection capability of Algorithm \ref{algorithm1}, we construct a privacy attacker that attempts to infer the original training images from the received information.
This attacker employs the iDLG attack model proposed in \cite{Zhao_iDLG}, which is more powerful than the well-known DLG attack model in \cite{Zhu_Deep} and can accurately reconstruct the original data from the shared gradients or model updates.
Fig. \ref{figure_idlg} shows that for the commonly used decentralized parallel stochastic gradient descent (D-PSGD) algorithm \cite{Lian_Can} in machine learning, the attacker can effectively recover the original training images from the shared model updates.
However, under the proposed algorithm with quantization, the attacker fails to infer any recognizable patterns, and the reconstructed images bear no resemblance to the original training samples.
This is because the quantization operation compresses the shared messages into a discrete set of values, irreversibly losing the precise information that is critical for reconstruction attacks like DLG or iDLG. Moreover, the quantization error accumulates over iterations, so each round of communication introduces additional distortion to the model updates. As a result, even after observing many rounds of communication, the attacker remains unable to recover the original images.
These results confirm the privacy protection capability of Algorithm \ref{algorithm1}, demonstrating that it can effectively preserve the local dataset of each agent.

\subsection{EEG-based Motor Imagery Classification on the ``BNCI2014001'' Dataset}
\begin{figure}[htp]
  \centering
  \includegraphics[width=0.45\textwidth]{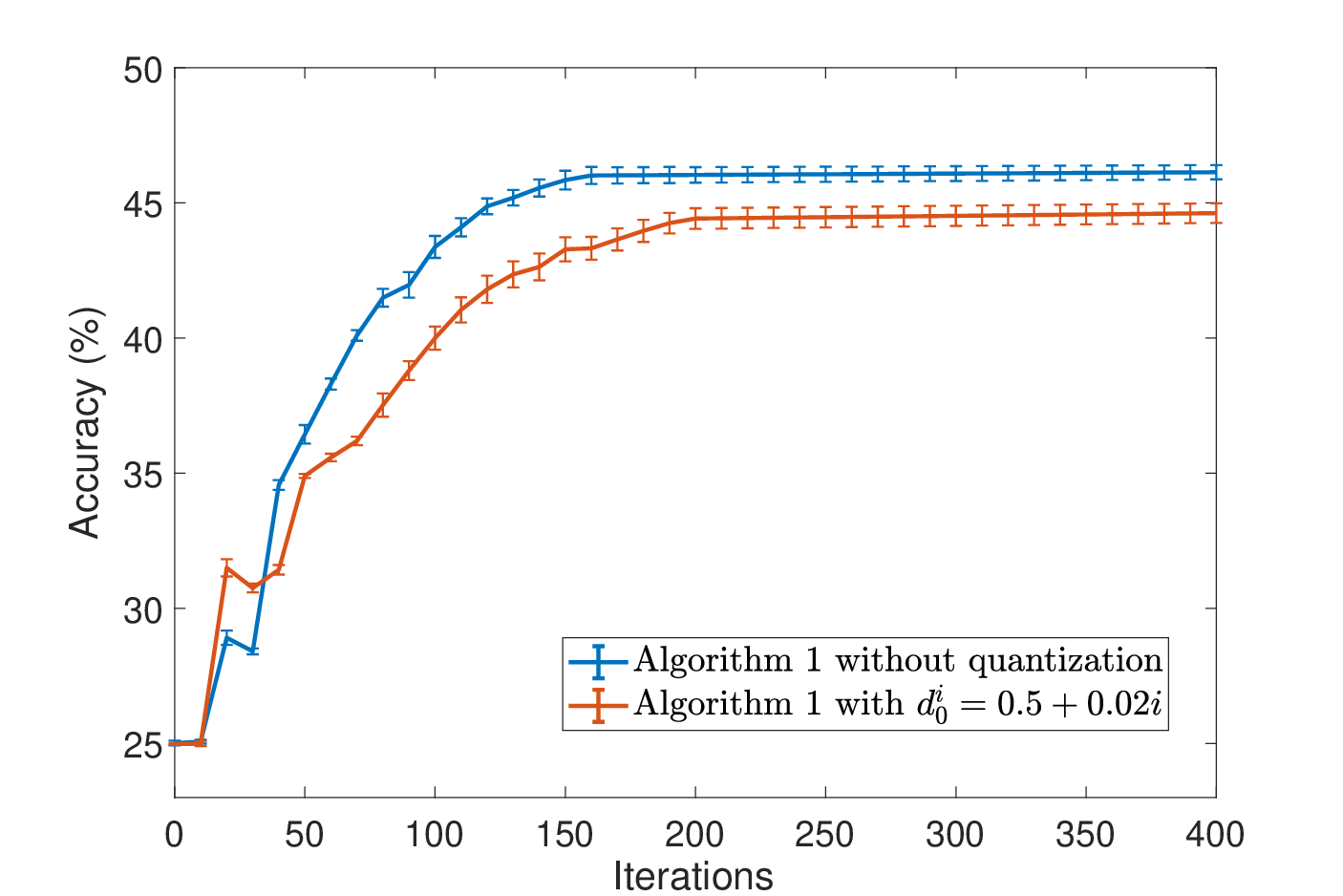}
  \caption{Comparison of the cross-subject classification accuracy of Algorithm 1 with and without quantization}
  \label{figure_EEG_acc}
\end{figure}

In the final experiment, we conduct motor imagery classification based on electroencephalography (EEG) via distributed online training of the EEG-based neural network (EEGNet) \cite{Lawhern_EEGNet} on the “BNCI2014001” dataset \cite{Tangermann_Review}. This dataset comprises four classes: left hand, right hand, both feet, and tongue, each with 144 trials \cite{Meng_User}. These EEG signals were recorded from 9 healthy subjects at 250 Hz with 22 channels, then bandpass filtered between 8 and 30 Hz with trials extracted over the [0, 4] s interval.
In the training, each agent maintains a local copy of the EEGNet whose hyperparameters are set the same as in \cite{Jia_Fderated}.
The performance evaluation follows a leave-one-subject-out cross-validation scheme, where one subject is designated as the test set, and the data from all remaining subjects are used for training.
Each subject was first performed with Euclidean alignment \cite{He_Transfer}.
In our distributed online training, each agent then randomly selects 20 trials from one or two fixed training subjects per iteration, with all training subjects exhaustively and disjointly assigned across agents.
To evaluate the proposed algorithm, the quantization stepsize and iteration stepsize are set as $d_t^i = \frac{0.5+0.02i}{(t+1)^{0.5+0.01i}}$ and $\lambda_t = \frac{0.2}{(t+1)^{0.61}}$ for $i=1,2,\cdots,5$, respectively.
The evolution of the cross-subject classification accuracies over 400 iterations for our algorithm and the non-privacy version is illustrated by the solid orange curve with error bars and the solid blue curve with error bars in Fig. \ref{figure_EEG_acc}, respectively. 
It can be observed that Algorithm \ref{algorithm1} with quantization converges to an accuracy close to the non-privacy version, which validates the resilience to quantization error and good generalization.
\begin{figure}[htp]
  \centering
  \includegraphics[width=0.47\textwidth]{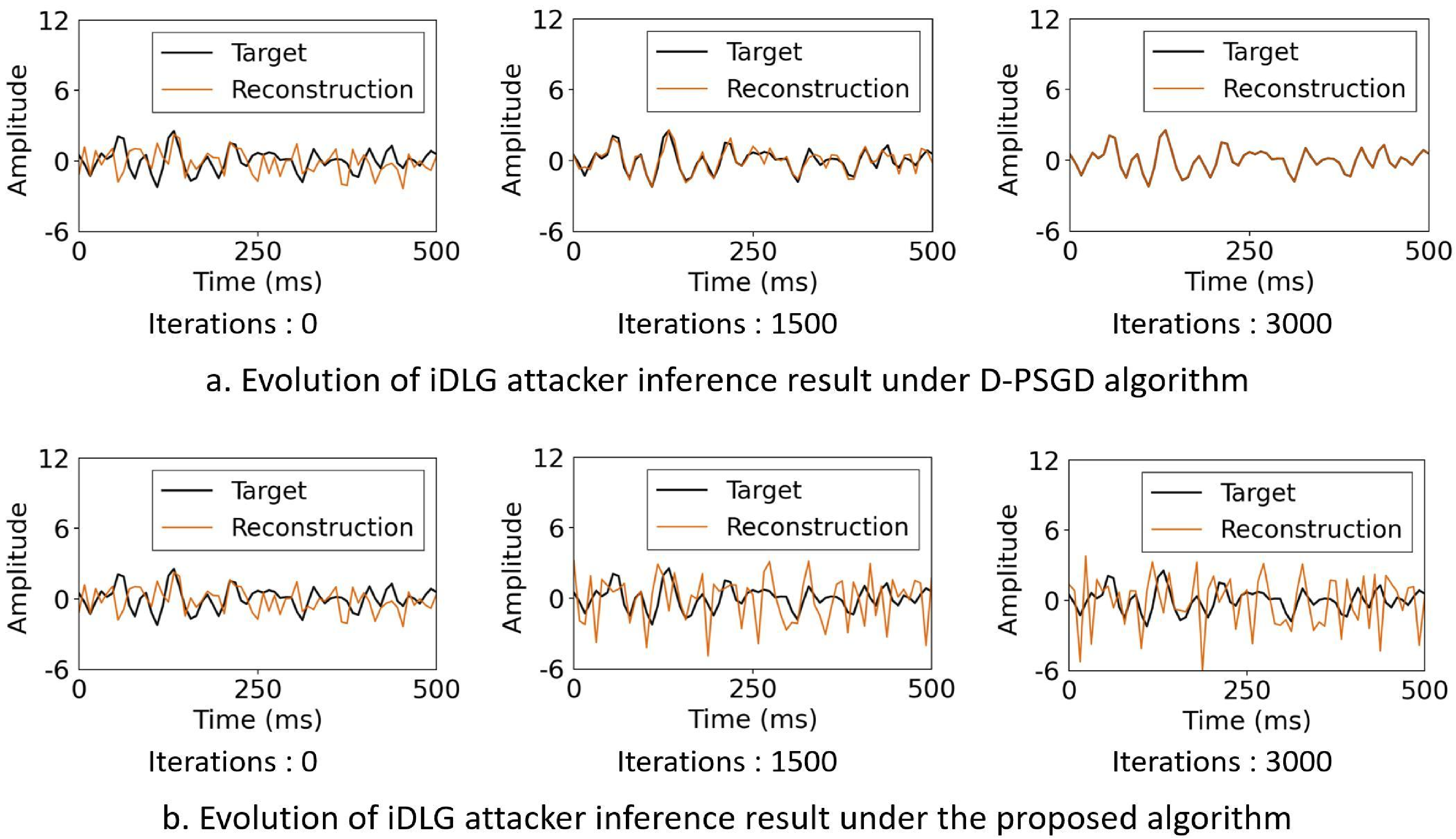}
  \caption{Comparison of iDLG attacker's EEG inference results under D-PSGD algorithm and our algorithm}
  \label{figure_idlg_EEG}
\end{figure}
\\
\indent\setlength{\parindent}{1em}
To verify the privacy protection capability of Algorithm
\ref{algorithm1} in the EEG-based motor imagery classification task, we again use the iDLG attack to restore the EEG data.
The result is shown in Fig. \ref{figure_idlg_EEG}, which displays only a portion of the reconstructed EEG data, specifically a single channel over a 0.5-second window.
It can be seen that without a privacy protection mechanism, EEG data for training by using D-PSGD can be effectively restored, while our algorithm is difficult to restore.
Furthermore, an attacker with access to identity-labeled EEG data can train a subject classifier that achieves a high recognition rate when matching the subject against reconstructed EEG data, due to the inherent correlation in the data.
Based on our experiments, the proposed algorithm is able to reduce this recognition rate from 74.3\% to 14.5\%. These results once again confirm the privacy protection capability of Algorithm \ref{algorithm1}, indicating that it can effectively protect each agent's local dataset and thereby successfully protect the subject's identity information in the EEG-based motor imagery classification task.

\section{Conclusion}
In this paper, we propose a quantization-based locally differentially private distributed online optimization algorithm.
In this algorithm, a stochastic quantizer with dynamic stepsize attenuation is designed to quantize the exchanged information of each agent prior to transmission, thereby protecting sensitive information.
Simultaneously, the online eigenvector estimation technique is utilized to enable our algorithm to operate in a fully distributed environment.
Theoretical analysis indicates that the proposed algorithm achieves almost sure convergence to the optimal solution while satisfying $(0,\delta^i)$-LDP, even as the number of iterations approaches infinity.
Finally, through numerical experimental results of distributed online training on the ``Mushrooms", ``MNIST", and ``BNCI2014001'' datasets, the effectiveness of the proposed method is verified.

\section*{Appendix}
The proof of Theorem 1 can be divided into four steps. In Step I-III, we analyze $\|\psi_t-\text{diag}(\boldsymbol{v})\boldsymbol{\bar{\psi}_t}\|_C^2$, $\|\theta_t-\boldsymbol{\bar{\theta}_t}\|_R^2$ and $F(\bar{\theta}_t)-F(\theta^*)$, respectively. In Step IV, we integrate these analyses to conclude the proof.
\\
\indent\setlength{\parindent}{1em}
\emph{Step I:} Analyze $\|\psi_t-\text{diag}(\boldsymbol{v})\boldsymbol{\bar{\psi}_t}\|_C^2$.
\\
\indent\setlength{\parindent}{1em}
By defining $\psi_t(\ell)=([\psi_t^1]_\ell, \cdots, [\psi_t^m]_\ell)^T$, $\theta_t(\ell)=([\theta_t^1]_\ell, \cdots, [\theta_t^m]_\ell)^T$, $\xi_{t,C}(\ell)=([\xi_{t,C}^1]_\ell, \cdots, [\xi_{t,C}^m]_\ell)^T$, $\zeta_{t,R}(\ell)=([\zeta_{t,R}^1]_\ell,\cdots,[\zeta_{t,R}^m]_\ell)^T$, and $g_t(\ell)=([g_t^1]_\ell, \cdots$,
$[g_t^m]_\ell)^T$, \eqref{vector_form} can be rewritten as the following per-coordinate form
\begin{align}
  \psi_{t+1}(\ell) &= \tilde{C}\psi_t(\ell) + \xi_{t,C}(\ell) + \lambda_t g_t(\ell), \nonumber\\
  \theta_{t+1}(\ell) &= \tilde{R} \theta_t + \zeta_{t,R}(\ell) - \tilde{z}_{t,m}^{-1}(\psi_{t+1}(\ell) - \psi_t(\ell)), \nonumber\\
  z_{t+1}(\ell) &= \tilde{R}z_t(\ell), \label{vector_form_ell}
\end{align}
for all $\ell=1,\cdots,d$, $t\ge 0$, where $\tilde{C}=I+C$, $\tilde{R}=I+R$ and $\tilde{z}_{t,m}=\text{diag}\{m[z_t^1]_1,\cdots,m[z_t^m]_m\}$. From the  definitions of $\bar{\psi}_t$ in \eqref{average_state_1} and $\psi_t(\ell)$ in \eqref{vector_form_ell}, we obtain $\psi_{t+1}(\ell)-v[\bar{\psi}_{t+1}]_\ell = \bar{C}(\psi_t(\ell) -v[\bar{\psi}_t]_\ell) + \Pi_v\xi_{t,C}(\ell) + \Pi_v\lambda_tg_t(\ell)$, where $\bar{C}=I+C-\frac{v\boldsymbol{1}^T}{m}$, $\Pi_v=I-\frac{v\boldsymbol{1}^T}{m}$, and the equation is obtained from $\bar{C}v=(I+C-\frac{v\boldsymbol{1}^T}{m})v=\boldsymbol{0}$.
\\
\indent\setlength{\parindent}{1em}
Taking the norm $\|\cdot\|_C^2$ on both sides yields
\begin{align}
  &\|\psi_{t+1}(\ell)-v[\bar{\psi}_{t+1}]_\ell\|_C^2 \nonumber\\
  &\leq (\|\bar{C}\|_C\|\psi_t(\ell)-v[\bar{\psi}_t]_\ell\|_C + \|\Pi_v\|_C \lambda_t \|g_t(\ell)\|_C)^2 \nonumber\\
  &\quad + 2\langle \bar{C}(\psi_t(\ell)-v[\bar{\psi}_t]_\ell) + \Pi_v \lambda_t g_t(\ell),\Pi_v \xi_{t,C}(\ell) \rangle_C \nonumber\\
  &\quad + \|\Pi_v\|_C^2 \|\xi_{t,C}(\ell)\|_C^2, \nonumber
\end{align}
where $\langle \cdot \rangle_C$ denotes the inner product induced by the norm $\| \cdot \|_C$.
Using the property $\|\bar{C}\|_C=1-\rho_C$ and for $\mathcal{F}_t = \{\theta_p^i,\psi_p^i,z_p^i,i\in[m],p\in[0,t]\}$, we take the conditional expectation $\mathbb{E}[\cdot|\mathcal{F}_t]$ on both sides and sum the relations over $\ell=1,\cdots,d$ to obtain
\begin{align}
  &\mathbb{E}[\|\psi_{t+1}-\text{diag}(\boldsymbol{v}) \boldsymbol{\bar{\psi}_{t+1}}\|_C^2 |\mathcal{F}_t] \nonumber\\
  &\quad\leq (1-\rho_C)\|\psi_t-\text{diag}(\boldsymbol{v})\boldsymbol{\bar{\psi}_t}\|_C^2 \nonumber\\
  &\qquad+ 2\|\Pi_v\|_C^2 \delta_{C,2}^2 \mathbb{E}[\|\xi_{t,C}\|^2|\mathcal{F}_t] \nonumber\\
  &\qquad + \Big(\frac{1}{\rho_C}+1\Big)\|\Pi_v\|_C^2\delta_{C,2} \lambda_t^2 \mathbb{E}[\|g_t\|^2|\mathcal{F}_t], \label{step1_sum}
\end{align}
in which the relations $2\langle a,b\rangle_C \leq \|a\|_C^2 + \|b\|_C^2$ and $\|a\|_C \leq \delta_{C,2}\|a\|_2$ ($\delta_{C,2}$ is the equivalence constant between the two norms) are used.
\\
\indent\setlength{\parindent}{1em}
For the second term in \eqref{step1_sum}, it can be derived from the property of quantizer \eqref{quantizer} that
\begin{align*}
 & \mathbb{E}[\|\xi_{t,C}\|^2 |\mathcal{F}_t] = \mathbb{E}\left[\sum_{i=1}^{m} \|\xi_{t,C}^i\|^2\big|\mathcal{F}_t\right] \nonumber\\
  &\quad\leq \sum_{i=1}^{m} \sum_{j\in \mathbb{N}_{C,i}^\text{in}} C_{ij}^2 \max_{i\in[m]}\frac{\{(d_t^i)^2\}}{4} 
  = \frac{C_C}{4} \max_{i\in[m]}\{(d_t^i)^2\}, 
\end{align*}
where $C_C = \sum_{i=1}^{m}\sum_{j\in\mathbb{N}_{C,i}^\text{in}}C_{ij}^2$.
By using Jensen's inequality and Assumption \ref{assumption 3}, \eqref{step1_sum} is deduced as
\begin{align}
  &\mathbb{E}[\|\psi_{t+1}-\text{diag}(\boldsymbol{v})\boldsymbol{\bar{\psi}_{t+1}}\|_C^2|\mathcal{F}_t] \nonumber\\
  &\quad\leq (1-\rho_C)\|\psi_t-\text{diag}(\boldsymbol{v})\boldsymbol{\bar{\psi}_t}\|_C^2 \nonumber\\
  &\qquad +\Big(\frac{1}{\rho_C}+1\Big)md_l^2\|\psi_v\|_C^2\delta_{C,2}^2\lambda_t^2 \nonumber\\
  &\qquad +\frac{\|\Pi_v\|_C^2\delta_{C,2}^2C_C}{2}\max_{i\in[m]}\{(d_t^i)^2\}. \label{step1_conculsion}
\end{align}
\\
\indent\setlength{\parindent}{1em}
\emph{Step II:} Analyze $\|\theta_t-\boldsymbol{\bar{\theta}_t}\|_R$.
\\
\indent\setlength{\parindent}{1em}
From the definitions of $\theta_t(\ell)$ in \eqref{vector_form_ell} and $\bar{\theta}_t$ in \eqref{average_state_1}, it is straightforward that
\begin{align}
  &\theta_{t+1}(\ell)-[\bar{\theta}_{t+1}]_\ell\boldsymbol{1} \nonumber\\
  &\quad= \tilde{R}\theta_t(\ell)-[\bar{\theta}_t]_\ell\boldsymbol{1} +\Big(I+\frac{\boldsymbol{1}u^T}{m}\Big)\zeta_{t,R}(\ell) \nonumber\\
  &\qquad -\Big(\tilde{U}^{-1}-\frac{\boldsymbol{1}\boldsymbol{1}^T}{m}\Big)(\psi_{t+1}(\ell)-\psi_t(\ell)) \nonumber\\
  &\qquad -\Big(I-\frac{\boldsymbol{1}u^T}{m}\Big)(\tilde{z}_{t,m}^{-1}-\tilde{U}^{-1})(\psi_{t+1}(\ell)-\psi_t(\ell)). \nonumber
\end{align}
Defining $\tilde{U}=\text{diag}(u)$, $\Pi_u=I-\frac{\boldsymbol{1}u^T}{m}$, $\Pi_U^e=(I-\frac{\boldsymbol{1}u^T}{m})(\tilde{z}_{t,m}^{-1}-\tilde{U}^{-1})$, $\Pi_U=\tilde{U}^{-1}-\frac{\boldsymbol{1}u^T}{m}$, $\Pi_z=\Pi_U^e+\Pi_U$, and using the property $\bar{R}(\theta_t(\ell)-[\bar{\theta}_t]_\ell\boldsymbol{1})=\tilde{R}\theta_t(\ell)-[\bar{\theta}_t]_\ell\boldsymbol{1}$, we arrive at $\theta_{t+1}(\ell)-[\bar{\theta}_{t+1}]_\ell\boldsymbol{1} = \bar{R}(\theta_t{\ell} -[\bar{\theta_t}]_\ell\boldsymbol{1}) +\Pi_u\zeta_{t,R}(\ell) 
  -\Pi_z(\psi_{t+1}(\ell)-\psi_t(\ell))$.
Using $Cv=\boldsymbol{0}$, it follows from \eqref{vector_form_ell} that
$\psi_{t+1}(\ell)-\psi_t(\ell) = C(\psi_t(\ell)-v[\bar{\psi}_t]_\ell)+\xi_{t,C}(\ell)+\lambda_tg_t(\ell). $
Then, we have
\begin{align}
  \theta_{t+1}(\ell)-[\bar{\theta}_t]_\ell\boldsymbol{1}&=\bar{R}(\theta_t(\ell)-[\bar{\theta}_t]_\ell\boldsymbol{1})  -\lambda_t\Pi_zg_t(\ell)\nonumber\\
  &\quad -\Pi_zC(\psi_t(\ell)-v[\bar{\psi}_t]_\ell)\nonumber\\
  &\quad+\Pi_u\zeta_{t,R}(\ell) - \pi_z\zeta_{t,C}(\ell). 
\end{align}
Taking $\|\cdot\|_R^2$ on both sides and using the inequalities $\|a-b\|_R \leq \|a\|_R +\|b\|_R$, $(a+b)^2\leq(1+\epsilon)a^2+\left(1+\frac{1}{\epsilon}\right)b^2$ with $\epsilon=\frac{\rho_R}{1-\rho_R}$, and the property $\|a\|_R \leq \delta_{R,2}\|a\|_2$, we take the conditional expectation $\mathbb{E}[\cdot|\mathcal{F}_t]$ on both sides. Then, summing the relations over $\ell=1,\cdots,d$ yields
\begin{align}
  &\mathbb{E}[\|\theta_{t+1}-\boldsymbol{\bar{\theta}_{t+1}}\|_R^2|\mathcal{F}_t] \leq (1-\rho_R)\|\theta_t-\boldsymbol{\bar{\theta}_t}\|_R^2\nonumber\\
  &\ \ +\frac{2\|\Pi_zC\|_R^2\delta_{R,C}^2}{\rho_R}\|\psi_t-\text{diag}(\boldsymbol{v})\boldsymbol{\bar{\psi}_t}\|_C^2 \nonumber\\
  &\ \ +4\delta_{R,2}^2(\|\Pi_u\|_R^2\mathbb{E}[\|\zeta_{t,R}\|^2|\mathcal{F}_t] +\|\Pi_z\|_R^2\mathbb{E}[\|\xi_{t,C}\|^2|\mathcal{F}_t]) \nonumber\\
  &\ \ +\Big(\frac{2}{\rho_R}+1\Big)\|\Pi_z\|_R^2\delta_{R,2}^2\lambda_t^2\mathbb{E}[\|g_t\|^2|\mathcal{F}_t]. \label{step2_E_vector}
\end{align}
For the second last term in \eqref{step2_E_vector}, it can be derived from the property of quantizer \eqref{quantizer} that $\mathbb{E}[\|\zeta_{t,R}^2|\mathcal{F}_t] \leq \frac{C_R}{4}\max_{i\in[m]}\{(d_t^i)^2\}$, where $C_R=\sum_{i=1}^{m}\sum_{j\in\mathbb{N}_{R,i}^\text{in}} R_{ij}^2$.
Substituting the result above into \eqref{step2_E_vector} yields
\begin{align}
  &\mathbb{E}[\|\theta_{t+1}-\boldsymbol{\bar{\theta}_{t+1}}\|_R^2|\mathcal{F}_t] \nonumber\\
  &\leq (1-\rho_R)\|\theta_t-\boldsymbol{\bar{\theta}_t}\|_R^2  +\left(\frac{2}{\rho_R}+1\right)md_l^2\|\Pi_z\|_R^2\delta_{R,2}^2\lambda_t^2 \nonumber\\
  &\quad +\frac{4\|C\|_R^2\delta_{R,C}^2(\|\Pi_U\|_R^2 +\|\Pi_U^e\|_R^2)}{\rho_R} \|\psi_t-\text{diag}(\boldsymbol{v})\boldsymbol{\bar{\psi}_t}\|_C^2 \nonumber\\
  &\quad +(\|\Pi_u\|_R^2C_R +\|\Pi_z\|_R^2C_C)\delta_{R,2}^2 \max_{i\in[m]}\{(d_t^i)^2\}, \label{step2_conculsion}
\end{align}
where the inequality is obtained by the property $\|\Pi_z\|_R^2 =\|\Pi_U+\Pi_U^e\|_R^2 \leq 2\|\Pi_U\|_R^2 +2\|\Pi_U^e\|_R^2$.
\\
\indent\setlength{\parindent}{1em}
\emph{Step III:} Analyze $F(\bar{\theta}_t)-F(\theta^*)$.
\\
\indent\setlength{\parindent}{1em}
From Jensen's inequality and Assumption \ref{assumption 2}, the following inequality holds for any $\theta_1,\theta_2\in\mathbb{R}^d$: $\|\nabla f_i(\theta_1)-\nabla f_i(\theta_2)\| \leq L\|\theta_1-\theta_2\|$ . Since $F(\theta)=\frac{1}{m}\sum_{i=1}^{m}f_i(\theta)$, it follows directly that $\|\nabla F(\theta_1) -\nabla F(\theta_2)\| \leq L\|\theta_1 - \theta_2\|$. Consequently, for any $a,b\in\mathbb{R}^d$, one has $F(a) \leq F(b) +\langle\nabla F(a), a-b\rangle +\frac{L}{2}\|a-b\|^2$.
By setting $a=\bar{\theta}_{t+1}$, $b=\bar{\theta}_t$, it is not difficult to obtain $F(\bar{\theta}_{t+1}) \leq F(\bar{\theta}_t) +\langle\nabla F(\bar{\theta}_{t+1}), \bar{\theta}_{t+1}-\bar{\theta}_t\rangle+\frac{L}{2}\|\bar{\theta}_{t+1}-\bar{\theta}_t\|^2$. Then, it follows from \eqref{average_state_1} that
\begin{align}
  \bar{\theta}_{t+1}-\bar{\theta}_t &= \bar{\zeta}_{t,R}-\Pi_{u,d}^e(\psi_{t+1}-\psi_t)-(\bar{\psi}_{t+1}-\bar{\psi}_t)  \nonumber\\
  &\quad -(\bar{\psi}_{t+1}-\bar{\psi}_t) +\Pi_{u,d}^e(\psi_t-\text{diag}(\boldsymbol{v})\boldsymbol{\bar{\psi}_t}) \nonumber\\
  &\quad-\Pi_{u,d}^e\text{diag}(\boldsymbol{v})(\boldsymbol{\bar{\psi}_{t+1}} -\boldsymbol{\bar{\psi}_t}), 
\end{align}
where $\Pi_{u,d}^e \triangleq \frac{1}{m}(u^T\otimes I_d)(z_{t,m}^{-1}-U^{-1})$. The last term of right-hand side of the above equation can be estimated as
$\Pi_{u,d}^e\text{diag}(\boldsymbol{v})(\boldsymbol{\bar{\psi}_{t+1}} -\boldsymbol{\bar{\psi}_t}) 
  = \frac{1}{m}(u^T(\tilde{z}_{t,m}^{-1}-\tilde{U}^{-1})\text{diag}(v)\boldsymbol{1_m})\otimes (\bar{\psi}_{t+1}-\bar{\psi}_t) 
  = \frac{1}{m}(u^T(\tilde{z}_{t,m}^{-1}-\tilde{U}^{-1})v)(\bar{\psi}_{t+1}-\bar{\psi}_t)$. 
Defining $\Pi_{u,v}^e \triangleq \frac{1}{m}(u^T(\tilde{z}_{t,m}^{-1}-\tilde{U}^{-1})v)$, one yields $\bar{\theta}_{t+1}-\bar{\theta}_t = \bar{\zeta}_{t,R} -(1+\Pi_{u,v}^e)\bar{\xi}_{t,C} -(1+\Pi_{u,v}^e)\lambda_t\bar{g}_t -\Pi_{u,d}^e(\psi_{t+1}-\text{diag}(\boldsymbol{v})\boldsymbol{\bar{\psi}_{t+1}}) +\Pi_{u,d}^e(\psi_t-\text{diag}(\boldsymbol{v})\boldsymbol{\bar{\psi}_t})$.
It is straightforward that
\begin{align}
  &\mathbb{E}[F(\bar{\theta}_{t+1})-F(\theta^*)|\mathcal{F}_t] \nonumber\\
  &\quad\leq F(\bar{\theta}_t)-F(\theta^*) -(1+\Pi_{u,v}^e)\lambda_t\mathbb{E}[\langle\nabla F(\bar{\theta}_t),\bar{g}_t\rangle|\mathcal{F}_t] \nonumber\\
  &\qquad -\mathbb{E}[\langle\nabla F(\bar{\theta}_t),\Pi_{u,d}^e(\psi_{t+1}-\text{diag}(\boldsymbol{v})\boldsymbol{\bar{\psi}_{t+1}})\rangle|\mathcal{F}_t] \nonumber\\
  &\qquad  +\frac{L}{2}\mathbb{E}[\|\bar{\zeta}_{t,R} -(1+\Pi_{u,v}^e)\bar{\xi}_{t,C} -(1+\Pi_{u,v}^e)\lambda_t\bar{g}_t  \nonumber\\
  &\qquad -\Pi_{u,d}^e(\psi_{t+1}-\text{diag}(\boldsymbol{v})\boldsymbol{\bar{\psi}_{t+1}})\|^2\nonumber\\
  &\qquad+\Pi_{u,d}^e(\psi_t-\text{diag}(\boldsymbol{v})\boldsymbol{\bar{\psi}_t})  |\mathcal{F}_t]\nonumber\\
  &\qquad+\langle\nabla F(\bar{\theta}_t),\Pi_{u,d}^e(\psi_t-\text{diag}(\boldsymbol{v})\boldsymbol{\bar{\psi}_t})\rangle. \label{step3_E}
\end{align}
Since $-\langle a,b\rangle =\frac{\|a-b\|^2-\|a\|^2-\|b\|^2}{2}$, it is derived that
\begin{align}
  &-\mathbb{E}[\langle\nabla F(\bar{\theta}_t),\bar{g}_t\rangle|\mathcal{F}_t] \nonumber\\
  &\quad= \frac{1}{2} \mathbb{E}[\|\nabla F(\bar{\theta}_t)-\bar{g}_t\|^2 -\|\nabla F(\bar{\theta}_t)\|^2 -\|\bar{g}_t\|^2|\mathcal{F}_t] \nonumber\\
  &\quad= \frac{1}{2} (\mathbb{E}[\|\nabla F(\bar{\theta}_t)-\bar{g}_t\|^2|\mathcal{F}_t] -\|\nabla F(\bar{\theta}_t)\|^2 -\|\bar{g}_t\|^2).\nonumber
\end{align}
Next, we deduce $\mathbb{E}[\|\nabla F(\bar{\theta}_t)-\bar{g}_t\|^2|\mathcal{F}_t]
\leq \frac{2}{m} \sum_{i=1}^{m}\mathbb{E}[$
$\|\nabla f_i(\theta_t^i) -\nabla f_t^i(\theta_t^i)\|^2|\mathcal{F}_t] +\frac{2}{m} \sum_{i=1}^{m}\|\nabla f_i(\bar{\theta}_t) -\nabla f_i(\theta_t^i)\|^2$.
By using Assumption \ref{assumption 2}, it can be directly inferred that $\mathbb{E}[\|\nabla f_i(\theta_t^i) -\nabla f_t^i(\theta_t^i)\|^2|\mathcal{F}_t]
= \frac{1}{(t+1)^2} \sum_{k=0}^{t} \mathbb{E}[ \|\nabla f_i(\theta_t^i) -\nabla l_i(\theta_t^i,x_k^i)\|^2 |\mathcal{F}_t] \leq \frac{\kappa^2}{t+1}$.
Therefore, it can be deduced that $\mathbb{E}[\|\nabla F(\bar{\theta}_t)-\bar{g}_t\|^2|\mathcal{F}_t] \leq \frac{2L^2}{m}\|\theta_t -\boldsymbol{\bar{\theta}_t}\|^2 +\frac{2\kappa^2}{t+1}$, which leads to 
\begin{align*}
  &-(1+\Pi_{u,v}^e)\lambda_t\mathbb{E}[\langle\nabla F(\bar{\theta}_t),\bar{g}_t\rangle|\mathcal{F}_t] \nonumber\\
  &\quad \leq \frac{L^2(1+\Pi_{u,v}^e)\delta_{2,R}^2}{m}\lambda_t\|\theta_t -\boldsymbol{\bar{\theta}_t}\|_R^2 \nonumber\\
  &\qquad +\frac{\kappa^2(1+\Pi_{u,v}^e)\lambda_t}{t+1} -\frac{(1+\Pi_{u,v}^e)\lambda_t}{2}\|\nabla F(\bar{\theta}_t)\|^2 \nonumber\\
  &\qquad -\frac{(1+\Pi_{u,v}^e)\lambda_t}{2}\|\bar{g}_t\|^2. \nonumber
\end{align*}
With the help of \eqref{step1_conculsion}, it can be obtained that
\begin{align}
  &-\mathbb{E}[\langle\nabla F(\bar{\theta}_t),\Pi_{u,d}^e(\psi_{t+1}-\text{diag}(\boldsymbol{v})\boldsymbol{\bar{\psi}_{t+1}})\rangle|\mathcal{F}_t] \nonumber\\
  &\quad\leq \frac{\|\Pi_c\|_C^2\delta_{C,2}^2\delta_{2,C}^2C_C}{4}\|\Pi_{u,d}^e\| \max_{i\in[m]}{(d_t^i)^2}\nonumber\\
  &\qquad +\frac{(1-\rho_C)\delta_{2,C}^2}{2}\|\Pi_{u,d}^e\|\|\psi_t -\text{diag}(\boldsymbol{v})\boldsymbol{\bar{\psi}_t}\|_C^2 \nonumber\\
  &\qquad +(\frac{1}{\rho_C}+1)\frac{m}{2}d_l^2\|\Pi_v\|_C^2\delta_{C,2}^2\delta_{2,C}^2\|\Pi_{u,d}^e\|\lambda_t^2\nonumber\\
  &\qquad+\frac{\|\Pi_{u,d}^e\|}{2} \|\nabla F(\bar{\theta}_t)\|^2.
\end{align}
Then, the following inequality holds
$\frac{L}{2}\mathbb{E}[\|\bar{\zeta}_{t,R} -(1+\Pi_{u,v}^e)\bar{\xi}_{t,C} +\Pi_{u,d}^e(\psi_t-\text{diag}(\boldsymbol{v})\boldsymbol{\bar{\psi}_t})-(1+\Pi_{u,v}^e)\lambda_t\bar{g}_t  -\Pi_{u,d}^e(\psi_{t+1}-\text{diag}(\boldsymbol{v})\boldsymbol{\bar{\psi}_{t+1}})\|^2 |\mathcal{F}_t]
  \leq \frac{5L}{2}(1+\Pi_{u,v}^e)^2\lambda_t^2\|\bar{g}_t\|^2$
  $+\frac{5L}{2}\|\Pi_{i,d}^e\|^2\|\psi_t -\text{diag}(\boldsymbol{v})\boldsymbol{\bar{\psi}_t}\|^2 
  +\frac{5L}{2}\mathbb{E}[\|\bar{\zeta}_{t,R}\|^2 |\mathcal{F}_t] +\frac{5L}{2}(1+\Pi_{u,v}^e)^2 \mathbb{E}[\|\bar{\xi}_{t,C}\|^2 |\mathcal{F}_t]
  +\frac{5L}{2}\|\Pi_{u,d}^e\|^2\mathbb{E}[\|\psi_{t+1} -\text{diag}(\boldsymbol{v})\boldsymbol{\bar{\psi}_{t+1}}\|^2 |\mathcal{F}_t]$.
By defining $\bar{C}_R \triangleq \frac{1}{m}\sum_{i=1}^{m} \sum_{j\in \mathbb{N}_{R,i}^\text{in}}$ $u_i^2 R_{ij}^2$ and $\bar{C}_C \triangleq \frac{1}{m}\sum_{i=1}^{m} \sum_{j\in \mathbb{N}_{C,i}^\text{in}} C_{ij}^2$, the preceding inequality further deduces $\mathbb{E}[\|\bar{\zeta}_{t,R}\|^2 |\mathcal{F}_t]\leq \frac{\bar{C}_R}{4} \max_{i\in[m]}\{(d_t^i)^2\}$ and $\mathbb{E}[\|\bar{\xi}_{t,C}\|^2 |\mathcal{F}_t] \leq \frac{\bar{C}_C}{4} \max_{i\in[m]}$ $\{(d_t^i)^2\}$.
Combining with the above analysis, we yield
\begin{align}
  &\frac{L}{2}\mathbb{E}[\|\bar{\zeta}_{t,R} -(1+\Pi_{u,v}^e)\bar{\xi}_{t,C} +\Pi_{u,d}^e(\psi_t-\text{diag}(\boldsymbol{v})\boldsymbol{\bar{\psi}_t}) \nonumber \\
  &-(1+\Pi_{u,v}^e)\lambda_t\bar{g}_t  -\Pi_{u,d}^e(\psi_{t+1}-\text{diag}(\boldsymbol{v})\boldsymbol{\bar{\psi}_{t+1}})\|^2 |\mathcal{F}_t] \nonumber\\
  &\leq \frac{5L}{2}(2-\rho_C)\delta_{2,C}^2\|\Pi_{u,d}^e\|^2\|\psi_t -\text{diag}(\boldsymbol{v})\boldsymbol{\bar{\psi}_t}\|_C^2 \nonumber \\
  &\quad +\frac{5L}{2}\Big(1+\frac{1}{\rho_C}\Big)md_l^2\|\Pi_v\|_C^2\delta_{C,2}^2\delta_{2,C}^2\|\Pi_{u,d}^e\|^2\lambda_t^2 \nonumber \\
  &\quad +\frac{5L}{2}(1+\Pi_{u,v}^e)^2\lambda_t^2\|\bar{g}_t\|^2 +\frac{5L}{8}\Big(\bar{C}_R +\bar{C}_C(1+\Pi_{u,v}^e)^2 \nonumber\\
  &\quad +2\|\Pi_v\|_C^2\delta_{C,2}^2\delta_{2,C}^2C_C\|\Pi_{u,d}^e\|^2\Big)\max_{i\in[m]}\{(d_t^i)^2\}.
\end{align}
For the last term of \eqref{step3_E}, one has $\langle\nabla F(\bar{\theta}_t),\Pi_{u,d}^e(\psi_t-\text{diag}(\boldsymbol{v})\boldsymbol{\bar{\psi}_t})\rangle \leq \frac{\|\Pi_{u,d}^e\|}{2}\|\nabla F(\bar{\theta}_t)\|^2
+\frac{\delta_{2,C}^2\|\Pi_{u,d}^e\|}{2}\|\psi_t -\text{diag}(\boldsymbol{v})\boldsymbol{\bar{\psi}_t}\|_C^2$.
Based on the above inequalities, we obtain
\begin{align}
  &\mathbb{E}[F(\bar{\theta}_{t+1})-F(\theta^*)|\mathcal{F}_t] \nonumber\\
  &\leq F(\bar{\theta}_t) -F(\theta^*) +( \frac{1}{m}L^2\delta_{2,R}^2\lambda_t+A_{t,1})\|\theta_t -\boldsymbol{\bar{\theta}_t}\|_R^2\nonumber\\
  &\  \ +A_{t,2} \|\psi_t -\text{diag}(\boldsymbol{v})\boldsymbol{\bar{\psi}_t}\|_C^2+\frac{\kappa^2(1+|\Pi_{u,v}^e|)}{t+1}\lambda_t\nonumber\\
  &\ \ +\frac{(\|\Pi_{u,d}^e\| +(|\Pi_{u,v}^e|-1)\lambda_t)}{2} \|\nabla F(\bar{\theta}_t)\|^2 \nonumber\\
  &\ \ +\frac{(\|\Pi_{u,d}^e\| +(|\Pi_{u,v}^e|-1)\lambda_t +5L(1+\Pi_{u,v}^e)^2\lambda_t^2)}{2} \|\bar{g}_t\|^2 \nonumber\\
  &\ \ + \frac{(1+5L\|\Pi_{u,d}^e\|)}{2}md_l^2\|\Pi_v\|_C^2\delta_{C,2}^2\delta_{2,C}^2\|\Pi_{u,d}^e\|\lambda_t^2 \nonumber\\
  &\ \ +(\frac{(1+5L\|\Pi_{u,d}^e\|)(\frac{1}{\rho_c}+1)}{4}\|\Pi_v\|_C^2\delta_{C,2}^2\delta_{2,C}^2C_C\|\Pi_{u,d}^e\| \nonumber\\
  &\  \ +\frac{5L}{8}(\bar{C}_R+(1+|\Pi_{u,v}^e|)^2\bar{C}_C))\max_{i\in[m]}\{(d_t^i)^2\},  \label{step3_conculsion}
\end{align}
where $A_{t,1}= \frac{1}{m}L^2\delta_{2,R}^2\lambda_t|\Pi_{u,v}^e|$ and $A_{t,2}= (1+5L\|\Pi_{u,d}^e\|)\left(1-\frac{\rho_C}{2}\right)\delta_{2,C}^2\|\Pi_{u,d}^e\|$.
\\
\indent\setlength{\parindent}{1em}
\emph{Step IV:} Combine Steps I-III and prove the theorem.
\\
\indent\setlength{\parindent}{1em}
By defining $\boldsymbol{v_t} =[F(\bar{\theta}_t) -F(\theta^*), \|\theta_t -\boldsymbol{\bar{\theta}_t}\|_R^2, \|\psi_t -\text{diag}(\boldsymbol{v})\boldsymbol{\bar{\psi}_t}\|_C^2]^T$, $\boldsymbol{u_t} =[\|\nabla F(\bar{\theta}_t)\|^2,\|\bar{g}_t\|^2]^T$ and $\bar{\rho}_C \triangleq \big(\frac{1}{\rho_C}+1\big)$,
it follows from \eqref{step1_conculsion}, \eqref{step2_conculsion}, and \eqref{step3_conculsion} that
\begin{align}
  \mathbb{E}[\boldsymbol{v_{t+1}} |\mathcal{F}_t] &\leq (V_t + A_t)\boldsymbol{v_t} -H_t\boldsymbol{u_t} +B_t, \nonumber\\
  &\leq (V_t + a_t \mathbf{11^T}) v_t + b_t \mathbf{1} - H_t \boldsymbol{u_t}, \label{lemma_2_form}
\end{align}
where $a_t \triangleq \max_{i,j\in[m]}\{[A_t]_{i,j}\}$, $b_t \triangleq \max_{i,j\in[m]}\{[B_t]_{i,j}\}$, and
\begin{align*}
  V_t &=
  \begin{bmatrix}
    1 & \frac{1}{m}L^2\delta_{2,R}^2\lambda_t & 0 \\
    0 & 1-\rho_R & \frac{4\|\Pi_UC\|_R^2\delta_{R,C}^2}{\rho_R} \\
    0 & 0 & 1-\rho_C
  \end{bmatrix}, \ B_t =
  \begin{bmatrix}
    B_{t,1} \\
    B_{t,2} \\
    B_{t,3}
  \end{bmatrix},\\
  A_t &=
  \begin{bmatrix}
    0 & A_{t,1} & A_{t,2} \\
    0 & 0 & A_{t,3} \\
    0 & 0 & 0
  \end{bmatrix}
  , \ H_t =
  \begin{bmatrix}
    H_{t,1} & H_{t,2} \\
    0 & 0 \\
    0 & 0
  \end{bmatrix},
\end{align*}
with $A_{t,3} = \frac{4}{\rho_R}\|C\|_R^2\delta_{R,C}^2\|\Pi_U^e\|_R^2$, 
$B_{t,1} = \kappa^2(1+|\Pi_{u,v}^e|)\frac{\lambda_t}{t+1} +\frac{(1+5L\|\Pi_{u,d}^e\|)}{2}\bar{\rho}_Cmd_l^2\|\Pi_v\|_C^2\delta_{C,2}^2\delta_{2,C}^2\|\Pi_{u,d}^e\|\lambda_t^2+(\frac{(1+5L\|\Pi_{u,d}^e\|)}{4}\|\Pi_v\|_C^2\delta_{C,2}^2\delta_{2,C}^2C_C\|\Pi_{u,d}^e\|+\frac{5L}{8}(\bar{C}_R +(1+|\Pi_{u,v}^e|)^2\bar{C}_C)\max_{i\in[m]}\{(d_t^i)^2\}$, 
$H_{t,1}= \frac{1}{2}\lambda_t -\frac{|\Pi_{u,v}^e|}{2}\lambda_t -\frac{\|\Pi_{u,d}^e\|}{2}$, 
$B_{t,2} = (\frac{2}{\rho_R}+1)md_l^2\|\Pi_z\|_R^2\delta_{R,2}^2\lambda_t^2+(\|\Pi_u\|_R^2C_R +\|\Pi_z\|_R^2C_C)\delta_{R,2}^2 \max_{i\in[m]}\{(d_t^i)^2\}$, $H_{t,2} = \frac{1}{2}\lambda_t -\frac{|\Pi_{u,v}^e|}{2}\lambda_t -\frac{\|\Pi_{u,d}^e\|}{2} -\frac{5L(1+|\Pi_{u,v}^e|)^2}{2}\lambda_t^2$, and $B_{t,3}= \bar{\rho}_Cmd_l^2\|\Pi_v\|_C^2\delta_{C,2}^2\lambda_t^2 +\frac{\|\Pi_v\|_C^2\delta_{C,2}^2C_C}{2} \max_{i\in[m]}\{(d_t^i)^2\}$.
\\
\indent\setlength{\parindent}{1em}
By the fact that the results of Lemma \ref{lemma 2} are asymptotic, they can be used when the starting index is shifted from $t=0$ to $t=T$, for any $T \ge 0$. Since $|\Pi_{u,v}^e|$ , $\|\Pi_{u,d}^e\|$ and $\lambda_t^2$ decay faster than $\lambda_t$, it can be easily found that all sequences in \eqref{lemma_2_form} are nonnegative for all $t \ge T$ (for some large enough $T \ge 0$).
\\
\indent\setlength{\parindent}{1em}
Similarly, it follows easily from Lemma \ref{lemma 4} that $|\Pi_{u,v}^e|$, $\|\Pi_{u,d}^e\|$, $\|\Pi_U^e\|_R^2$ are all summable, while $\frac{\lambda_t}{t+1}$, $\lambda_t^2$ and $(d_t^i)^2$ are summable sequences. Consequently, ${a_t}$ and ${b_t}$ are both summable.
Then, we tend to prove that there exists a vector $\pi > 0$ such that $\pi^T V_t \leq \pi^T$ and $\pi^T H_t \ge 0 $ hold a.s. for all $t \ge 0$. To this end, it suffices to verify the following inequalities
\begin{align}
  &\frac{1}{m}L^2\delta_{2,R}^2\lambda_t\pi_1 +(1-\rho_R)\pi_2 \leq \pi_2, \nonumber\\
  &\frac{4}{\rho_R}\|\Pi_UC\|_R^2\delta_{R,C}^2\pi_2 +(1-\rho_C)\pi_3 \leq \pi_3. \label{lemma_pi}
\end{align}
\\
\indent\setlength{\parindent}{1em}
Note that the first inequality in \eqref{lemma_pi} equals to $\pi_2 \ge \frac{1}{m\rho_R}L^2\delta_{2,R}^2\lambda_t\pi_1$. Given that $\lim_{t\to\infty}\lambda_t=0$ holds and all parameters on the right side are positive, it can be easily obtained that for a given $\pi_1>0$, there is a $\pi_2>0$ satisfying the above inequality for $t \ge T$.
In \eqref{lemma_pi}, the second inequality is equal to $\pi_3 \ge \frac{4\|\Pi_uC\|_R^2\delta_{R,C}^2}{\rho_R\rho_C}\pi_2$. For all parameters on the right side are positive, it can be seen that for a given $\pi_2>0$, there is a $\pi_3>0$ satisfying the inequality for $t \ge T$.
In addition, Since $h_t(1)$ and $h_t(2)$ are both nonnegative for $t \ge T$, it is not difficult to obtain that $\pi^TH_t \ge 0$ holds for $t \ge T$. Thus, we can always find a vector $\pi$ satisfying the inequalities in \eqref{lemma_pi} for $t\ge T$ with some large enough $T\ge 0$. Hence, the conditions in Lemma \ref{lemma 2} are satisfied.
\\
\indent\setlength{\parindent}{1em}
Based on the above results, the two conclusions of Theorem \ref{theorem 1} can be proved.
\\
\indent\setlength{\parindent}{1em}
a) By Lemma \ref{lemma 2}, it can be obtained that $\lim_{t\to\infty}\pi^T\boldsymbol{v_t}$ exists a.s. and $\sum_{t=0}^{\infty}\pi^TH_t\boldsymbol{u_t}<\infty$ holds a.s.. Since $\pi>0$ and $\boldsymbol{v_t}\ge 0$, it follows that $\lim_{t\to\infty}F(\bar{\theta}_t)$ exists a.s.. Since $\pi^TH_t =[\pi_1h_t(1), \pi_1h_t(2)]$, and $\|\Pi_{u,v}^e\|$, $\|\Pi_{u,d}^e\|$, and $\lambda_t^2$ are summable, one has
\begin{equation}
  \sum_{t=0}^{\infty}\lambda_t\|\nabla F(\bar{\theta}_t)\|^2 <\infty,\ \sum_{t=0}^{\infty}\lambda_t\|\bar{g}_t\|^2 <\infty\ \ \text{a.s..} \label{nabla F}
\end{equation}
Now, to prove that both $\|\theta_t^i-\bar{\theta}_t\|$ and $\|\psi_t^i -v_i\bar{\psi}_t\|$ converge to $0$ a.s., from \eqref{step1_conculsion} one lets $v_t=\|\pi_t-\text{diag}(\boldsymbol{v})\boldsymbol{\bar{\psi}_t}\|_C^2$, $\alpha_t=0$, $q_t=\rho_C$, and the remainder as $p_t$. Hence, all conditions in Lemma \ref{lemma 1} are satisfied, and it follows that
\begin{equation}
  \sum_{t=0}^{\infty}\|\psi_t^i -v_i\bar{\psi}_t\|^2 <\infty,\ \lim_{t\to\infty}\|\psi_t^i -v_i\bar{\psi}_t\|=0\ \ \text{a.s..} \label{lim psi}
\end{equation}
Similarly, from \eqref{step2_conculsion} one lets $v_t=\|\theta_t-\boldsymbol{\bar{\theta}_t}\|_R^2$, $\alpha_t=0$, $q_t=\rho_R$, and the remainder as $p_t$. Hence, all conditions in Lemma \ref{lemma 1} are valid and the following is obtained
\begin{equation}
  \sum_{t=0}^{\infty}\|\theta_t^i -\bar{\theta}_t\|^2<\infty,\ \lim_{t\to\infty}\|\theta_t^i -\bar{\theta}_t\|=0\ \ {a.s..}  \label{lim theta}
\end{equation}
\\
\indent\setlength{\parindent}{1em}
b) From \eqref{nabla F}, since $\sum_{t=0}^{\infty}\lambda_t=\infty$, it follows directly that $\liminf_{t\to\infty}\|\nabla F(\bar{\theta}_t)\|=0$ a.s.. In addition, since $F(\cdot)$ has bounded level sets, then the sequence $\{\bar{\theta}_t\}$ is a.s. bounded and thus possesses accumulation points. Let $\{\bar{\theta}_t^i\}$ be a subsequence such that $\lim_{i\to\infty}\|\nabla F(\bar{\theta}_t^i)\|=0$ a.s.. Without loss of generality, we may assume that $\{\bar{\theta}_t^i\}$ converges a.s.. Otherwise, one can extract a convergent subsequence of $\{\bar{\theta}_t^i\}$. Denoting $\lim_{i\to\infty}\bar{\theta}_t^i=\hat{\theta}$, the continuity of $\nabla F(\cdot)$ implies $\nabla F(\hat{\theta})=0$, so $\hat{\theta}$ is an optimal point. By the continuity of $F(\cdot)$, it follows that $F(\hat{\theta})=F^*$. Since $\lim_{t\to\infty}F(\bar{\theta}_t)$ exists a.s., we conclude that $\lim_{t\to\infty}F(\bar{\theta}_t)=F^*$ a.s..
\\
\indent\setlength{\parindent}{1em}
Simultaneously, from \eqref{lim theta} one has $\lim_{t\to\infty}\|\theta_t^i-\bar{\theta}_t\|=0$ a.s. for all $i\in[m]$. Consequently, each $\{\theta_t^i\}$ shares the same accumulation points as $\{\bar{\theta}_t\}$ a.s., implying that $\lim_{t\to\infty}F(\theta_t^i)=F^*$ a.s. for all $i\in[m]$ according to the continuity of $F(\cdot)$. The proof of Theorem 1 is complete.

\end{document}